\documentclass[11pt,twoside]{article}

\usepackage{amsfonts,epsfig,graphicx,subfigure}
\usepackage{amsmath,amssymb,amsthm}

\usepackage{epic,eepic}

\usepackage{epsf}
\usepackage{fancyheadings}
\usepackage{graphics}
\usepackage{amsfonts}
\usepackage{amsmath}
\usepackage{amssymb}
\usepackage{epic}
\usepackage{eepic}
\usepackage{eepicemu}

\newenvironment{carlist}
 {\begin{list}{$\bullet$}
 {\setlength{\topsep}{0in} \setlength{\partopsep}{0in}
  \setlength{\parsep}{0in} \setlength{\itemsep}{\parskip}
  \setlength{\leftmargin}{0.07in} \setlength{\rightmargin}{0.08in}
  \setlength{\listparindent}{0in} \setlength{\labelwidth}{0.08in}
  \setlength{\labelsep}{0.1in} \setlength{\itemindent}{0in}}}
 {\end{list}}

\newcommand{\bcar}{\begin{carlist}}
\newcommand{\ecar}{\end{carlist}}

\newenvironment{carliste}
 {\begin{list}x
 {\setlength{\topsep}{0in} \setlength{\partopsep}{0in}
  \setlength{\parsep}{0in} \setlength{\itemsep}{\parskip}
  \setlength{\leftmargin}{0.07in} \setlength{\rightmargin}{0.08in}
  \setlength{\listparindent}{0in} \setlength{\labelwidth}{0.08in}
  \setlength{\labelsep}{0.1in} \setlength{\itemindent}{0in}}}
 {\end{list}}

\newcommand{\bcare}{\begin{carliste}}
\newcommand{\ecare}{\end{carliste}}

\newcommand{\Prob}{\ensuremath{\mathbb{P}}}

\makeatletter
\long\def\@makecaption#1#2{
        \vskip 0.8ex
        \setbox\@tempboxa\hbox{\small {\bf #1:} #2}
        \parindent 1.5em  
        \dimen0=\hsize
        \advance\dimen0 by -3em
        \ifdim \wd\@tempboxa >\dimen0
                \hbox to \hsize{
                        \parindent 0em
                        \hfil 
                        \parbox{\dimen0}{\def\baselinestretch{0.96}\small
                                {\bf #1.} #2
                                } 
                        \hfil}
        \else \hbox to \hsize{\hfil \box\@tempboxa \hfil}
        \fi
        }
\makeatother

\long\def\comment#1{}

\makeatletter
\def\@cite#1#2{[\if@tempswa #2 \fi #1]}
\makeatother

\makeatletter
\long\def\barenote#1{
    \insert\footins{\footnotesize
    \interlinepenalty\interfootnotelinepenalty 
    \splittopskip\footnotesep
    \splitmaxdepth \dp\strutbox \floatingpenalty \@MM
    \hsize\columnwidth \@parboxrestore
    {\rule{\z@}{\footnotesep}\ignorespaces
      #1\strut}}}
\makeatother

\newcommand{\bit}{\begin{itemize}}
\newcommand{\eit}{\end{itemize}}
\newcommand{\ben}{\begin{enumerate}}
\newcommand{\een}{\end{enumerate}}

\newcommand{\bear}{\begin{eqnarray}}
\newcommand{\eear}{\end{eqnarray}}

\newcommand{\vtiny}{\vspace*{.1in}}

\newcommand{\ones}{\ensuremath{\mathbf{1}}}

\newcommand{\Exs}{\ensuremath{{\mathbb{E}}}}

\newcommand{\beq}{\begin{quotation}}
\newcommand{\enq}{\end{quotation}}

\newcommand{\estart}{\begin{equation}}
\newcommand{\eend}{\end{equation}}

\newcommand{\widgraph}[2]{\includegraphics[keepaspectratio,width=#1]{#2}}

\newcommand{\defn}{\ensuremath{:  =}}

\newcommand{\ysca}{{{y}}}

\newcommand{\wsca}{{{w}}}

\newcommand{\Ysca}{{{Y}}}

\newcommand{\Wsca}{{{W}}}

\newcommand{\bec}{\begin{center}}
\newcommand{\enc}{\end{center}}

\newcommand{\beit}{\begin{itemize}}
\newcommand{\enit}{\end{itemize}}

\newcommand{\been}{\begin{enumerate}}
\newcommand{\enen}{\end{enumerate}}

\newcommand{\comsl}{\begin{slide}}
\newcommand{\comspor}{\begin{slide*}}

\newcommand{\comsld}[2]{\begin{slide}[#1,#2]}
\newcommand{\comspord}[2]{\begin{slide*}[#1,#2]}

\newcommand{\mendsl}{\end{slide}}
\newcommand{\mendspo}{\end{slide*}}

\newcommand{\estim}[1]{\ensuremath{\widehat{#1}}}

\newcommand{\wtil}[1]{\ensuremath{\widetilde{#1}}}

\newcommand{\real}{\ensuremath{{\mathbb{R}}}}

\DeclareMathOperator{\var}{var}

\DeclareMathOperator{\cov}{cov}

\DeclareMathOperator{\trace}{trace}

\theoremstyle{plain}

\newtheorem{theo}{Theorem}[section]

\newtheorem{lem}{Lemma}[section]
\newtheorem{prop}{Proposition}[section]
\newtheorem{cor}{Corollary}[section]

\theoremstyle{definition} 

\newtheorem{nota}{Notation}[section]
\newtheorem{de}{Definition}[section]
\newtheorem{exa}{Example}[section]
\newtheorem{as}{Assumption}[section]
\newtheorem{alg}{Algorithm}[section]

\newcommand{\btheo}{\begin{theo}}
\newcommand{\bde}{\begin{de}}
\newcommand{\ble}{\begin{lem}}
\newcommand{\bpr}{\begin{prop}}
\newcommand{\bno}{\begin{nota}}
\newcommand{\bex}{\begin{exa}}
\newcommand{\bcor}{\begin{cor}}
\newcommand{\spro}{\begin{proof}}
\newcommand{\bas}{\begin{as}}
\newcommand{\balg}{\begin{alg}}

\newcommand{\etheo}{\end{theo}}
\newcommand{\ede}{\end{de}}
\newcommand{\ele}{\end{lem}}
\newcommand{\epr}{\end{prop}}
\newcommand{\eno}{\end{nota}}
\newcommand{\eex}{\end{exa}}
\newcommand{\ecor}{\end{cor}}
\newcommand{\fpro}{\end{proof}}
\newcommand{\eas}{\end{as}}
\newcommand{\ealg}{\end{alg}}

\theoremstyle{plain}

\newtheorem{theos}{Theorem}
\newtheorem{props}{Proposition}
\newtheorem{lems}{Lemma}
\newtheorem{cors}{Corollary}

\theoremstyle{definition}
\newtheorem{exas}{Example}
\newtheorem{algs}{Algorithm}
\newtheorem{asss}{Asumption}
\newtheorem{defns}{Definition}

\newcommand{\btheos}{\begin{theos}}
\newcommand{\etheos}{\end{theos}}
\newcommand{\bprops}{\begin{props}}
\newcommand{\eprops}{\end{props}}
\newcommand{\bdes}{\begin{defns}}
\newcommand{\edes}{\end{defns}}
\newcommand{\blems}{\begin{lems}}
\newcommand{\elems}{\end{lems}}
\newcommand{\bcors}{\begin{cors}}
\newcommand{\ecors}{\end{cors}}
\newcommand{\bexs}{\begin{exas}}
\newcommand{\eexs}{\end{exas}}
\newcommand{\balgs}{\begin{algs}}
\newcommand{\ealgs}{\end{algs}}
\newcommand{\bass}{\begin{asss}}
\newcommand{\eass}{\end{asss}}

\newcommand{\myeigmin}{\ensuremath{\Lambda_{min}}}
\newcommand{\myeigmax}{\ensuremath{\Lambda_{max}}}

\newcommand{\Sset}{\ensuremath{S}}

\newcommand{\Amat}{\ensuremath{X}}
\newcommand{\Amatt}[1]{\ensuremath{\Amat_{#1}}}
\newcommand{\arow}{\ensuremath{x}}
\newcommand{\acol}{\ensuremath{X}}

\newcommand{\betastar}{\ensuremath{\beta^*}}
\newcommand{\betasca}{\ensuremath{\beta}}
\newcommand{\betahat}{\ensuremath{\widehat{\beta}}}

\newcommand{\esca}{\ensuremath{w}}
\newcommand{\Esca}{\ensuremath{W}}

\newcommand{\SolSet}{\ensuremath{\mathcal{S}}}

\newcommand{\Sbar}{\ensuremath{S^c}}
\newcommand{\sgn}{\ensuremath{\operatorname{sgn}}}

\newcommand{\Prop}{\ensuremath{\mathcal{R}}}

\newcommand{\vstar}{\ensuremath{v^*}}

\newcommand{\range}{\ensuremath{\operatorname{range}}}

\newcommand{\subgrad}{\ensuremath{z}}
\newcommand{\subgradopt}{\ensuremath{\estim{\subgrad}}}

\newcommand{\svec}{\ensuremath{\spec}}

\newcommand{\sigw}{\ensuremath{\sigma^2}}
\newcommand{\sigwfour}{\ensuremath{\sigma^4}}

\newcommand{\Uvar}{\ensuremath{U}}
\newcommand{\Vvar}{\ensuremath{V}}
\newcommand{\Vvarcond}{\ensuremath{\wtil{\Vvar}}}

\newcommand{\Zvar}{\ensuremath{Z}}
\newcommand{\Zvarcond}{\ensuremath{\wtil{\Zvar}}}

\newcommand{\spec}{\ensuremath{\vec{b}\,}}
\newcommand{\myVar}{\ensuremath{M_\numobs}}
\newcommand{\myVarr}{\ensuremath{H}}
\newcommand{\myTemp}{\ensuremath{T}}

\newcommand{\evA}{\ensuremath{\mathcal{M}(\Vvar)}}
\newcommand{\evB}{\ensuremath{\mathcal{M}(\Uvar)}}

\newcommand{\Cmin}{\ensuremath{C_{min}}}
\newcommand{\Cmax}{\ensuremath{C_{max}}}

\newcommand{\CovMat}{\ensuremath{\Sigma}}
\newcommand{\CondMat}{\ensuremath{\CovMat_{(\Sbar  |  S)}}}

\newcommand{\ThreshLow}{\ensuremath{\theta_{\ell}}}
\newcommand{\ThreshUp}{\ensuremath{\theta_{u}}}

\newcommand{\Tail}{\ensuremath{\mathcal{T}}}

\newcommand{\interdelta}{\ensuremath{\delta'}}
\newcommand{\threshbou}{\ensuremath{\nu}}

\newcommand{\numobs}{\ensuremath{n}}
\newcommand{\spindex}{\ensuremath{s}}
\newcommand{\mdim}{\ensuremath{p}}

\newcommand{\bigN}{\ensuremath{N}}

\newcommand{\Vtil}{\ensuremath{\wtil{\Vvar}}}
\newcommand{\Vtilmax}{\ensuremath{\Vtil_{(\bigN)}}}

\newcommand{\Uvartil}{\ensuremath{\wtil{\Uvar}}}

\newcommand{\matbound}{\ensuremath{\epsilon}}

\newcommand{\toep}{\ensuremath{\rho}}

\newcommand{\myumean}{\ensuremath{Y}}

\newcommand{\myuvar}{\ensuremath{Y'}}

\newcommand{\mustar}{\ensuremath{\mu^*}}
\newcommand{\Dconmax}{\ensuremath{D_{\operatorname{max}}}}

\begin{document}
\begin{center}

{{\LARGE \bf{Sharp thresholds for high-dimensional and noisy recovery
of sparsity}}}

\vspace*{.5in}

{\large {
\begin{center}
Martin J. Wainwright \\
Department of Statistics, and \\
Department of Electrical Engineering and Computer Sciences \\
University of California, Berkeley \\
\texttt{wainwrig@\{eecs,stat\}.berkeley.edu}
\end{center}

}}

\vtiny

{\large Technical Report, UC Berkeley, Department of Statistics} \\
{\large May 2006}

\vtiny

\end{center}

\begin{abstract}
The problem of consistently estimating the sparsity pattern of a
vector $\betastar \in \real^\mdim$ based on observations contaminated
by noise arises in various contexts, including subset selection in
regression, structure estimation in graphical models, sparse
approximation, and signal denoising. We analyze the behavior of
$\ell_1$-constrained quadratic programming (QP), also referred to as
the Lasso, for recovering the sparsity pattern.  Our main result is to
establish a sharp relation between the problem dimension $\mdim$, the
number $\spindex$ of non-zero elements in $\betastar$, and the number
of observations $\numobs$ that are required for reliable recovery.
For a broad class of Gaussian ensembles satisfying mutual incoherence
conditions, we establish existence and compute explicit values of
thresholds $\ThreshLow$ and $\ThreshUp$ with the following properties:
for any $\threshbou > 0$, if $\numobs > 2 \, ( \ThreshUp + \threshbou)
\log (\mdim - \spindex) + \spindex + 1$, then the Lasso succeeds in
recovering the sparsity pattern with probability converging to one for
large problems, whereas for $\numobs < 2 \, ( \ThreshLow - \threshbou)
\log (\mdim - \spindex) + \spindex + 1$, then the probability of
successful recovery converges to zero.  For the special case of the
uniform Gaussian ensemble, we show that $\ThreshLow = \ThreshUp = 1$,
so that the threshold is sharp and exactly determined.
\end{abstract}


\section{Introduction}

The problem of recovering the sparsity pattern of an unknown vector
$\betastar$---that is, the positions of the non-zero entries of
$\betastar$--- based on noisy observations arises in a broad variety
of contexts, including subset selection in regression~\cite{Miller90},
structure estimation in graphical models~\cite{Meinshausen06}, sparse
approximation~\cite{Devore93,Natarajan95}, and signal
denoising~\cite{Chen98}.  A natural optimization-theoretic formulation
of this problem is via \mbox{$\ell_0$-minimization,} where the
$\ell_0$ ``norm'' of a vector corresponds to the number of non-zero
elements.  Unfortunately, however, $\ell_0$-minimization problems are
known to be NP-hard in general~\cite{Natarajan95}, so that the
existence of polynomial-time algorithms is highly unlikely.  This
challenge motivates the use of computationally tractable
approximations or relaxations to $\ell_0$ minimization.  In
particular, a great deal of research over the past decade has studied
the use of the $\ell_1$-norm as a computationally tractable surrogate
to the $\ell_0$-norm.

In more concrete terms, suppose that we wish to estimate an unknown
but fixed vector $\betastar \in \real^{\mdim}$ on the basis of a set
of $n$ observations of the form
\begin{eqnarray}
\label{EqnLinearObs}
\Ysca_k & = & \arow_k^T \betastar + \Wsca_k, \qquad k = 1, \ldots
\numobs,
\end{eqnarray}
where $\arow_k \in \real^{\mdim}$, and $\Wsca_k \sim N(0, \sigw)$ is
additive Gaussian noise.  In many settings, it is natural to assume
that the vector $\betastar$ is \emph{sparse}, in that its
\emph{support}
\begin{eqnarray}
\label{EqnSupport}
\Sset & \defn & \{ i \in \{1, \ldots \mdim \} \; \mid \; \betastar_i
\neq 0 \}
\end{eqnarray}
has relatively small cardinality $s = |\Sset|$.  Given the observation
model~\eqref{EqnLinearObs} and sparsity assumption~\eqref{EqnSupport},
a reasonable approach to estimating $\betastar$ is by solving the
$\ell_1$-constrained quadratic program (QP)
\begin{equation}
\label{EqnOrigQP}
\min_{\beta \in \real^\mdim} \left \{\frac{1}{2 \numobs} \sum_{k=1}^n
\| \Ysca_k - \arow^T_k \beta \|_2^2 + \lambda_{\numobs} \| \beta \|_1
\right \},
\end{equation}
where $\lambda_{\numobs} \geq 0$ is a regularization parameter.  Of interest
are conditions on the \emph{ambient dimension} $\mdim$, the
\emph{sparsity index} $s$, and the \emph{number of observations} $n$
for which it is possible (or impossible) to recover the support set
$\Sset$ of $\betastar$.

\subsection{Overview of previous work}

Given the substantial literature on the use of $\ell_1$ constraints
for sparsity recovery and subset selection, we provide only a very
brief (and hence necessarily incomplete) overview here.  In the
\emph{noiseless version} ($\sigma^2 = 0$) of the linear observation
model~\eqref{EqnLinearObs}, one can imagine estimating $\betastar$ by
solving the problem
\begin{equation}
\label{EqnLPRelax}
\min_{\beta \in \real^\mdim} \|\beta \|_1 \qquad \mbox{subject to}
\quad \arow^T_k \beta = \Ysca_k, \quad k = 1, \ldots, \numobs.
\end{equation}
This problem is in fact a linear program (in disguise), and
corresponds to a method in signal processing known as basis pursuit,
pioneered by Chen et al.~\cite{Chen98}.  For the noiseless setting,
the interesting regime is the underdetermined setting (i.e., $\numobs
< \mdim$).  With contributions from a broad range of
researchers~\cite[e.g.,]{CanRomTao04,Chen98,Donoho01,Donoho06,Elad02,Feuer03,Malioutov04,Tropp04},
there is now a fairly complete understanding of conditions on
deterministic vectors $\{\arow_k\}$ and sparsity index $\spindex$ for
which the true solution $\betastar$ can be recovered exactly.  Without
going into technical details, the rough idea is that the \emph{mutual
incoherence} of the vectors $\{\arow_k\}$ must be large relative to
the sparsity index $\spindex$, and indeed we impose similar conditions
to derive our results (e.g., conditions~\eqref{EqnKeyCondA}
and~\eqref{EqnCovCond} in the sequel).  Most closely related to the
current paper---as we discuss in more detail in the sequel---are
recent results by Donoho~\cite{Donoho04a}, as well as Candes and
Tao~\cite{CandesTao05} that provide high probability results for
random ensembles.  More specifically, as independently established by
both sets of authors using different methods, for uniform Gaussian
ensembles (i.e., $\arow_k \sim N(0, I_{\mdim})$) with the ambient
dimension $\mdim$ scaling linearly in terms of the number of
observations (i.e., $\mdim = \gamma \numobs$, for some $\gamma > 1$),
there exists a constant $\alpha > 0$ such that all sparsity patterns
with $s \leq \alpha \mdim$ can be recovered with high probability.

There is also a substantial body of work focusing on the noisy setting
($\sigma^2 > 0$), and the use of quadratic programming techniques for
sparsity
recovery~\cite[e.g.,]{Chen98,Fuchs04,Fuchs05,Tropp06,DonElTem06,Fletcher06,Meinshausen06,Zhao06}. The
$\ell_1$-constrained quadratic program~\eqref{EqnOrigQP}, also known
as the Lasso~\cite{Tibshirani96,Efron04}, has been the focus of
considerable research in recent years.  Knight and Fu~\cite{Knight00}
analyze the asymptotic behavior of the optimal solution, not only for
$\ell_1$ regularization but for $\ell_p$-regularization with $p \in
(0,2]$.  Fuchs~\cite{Fuchs04,Fuchs05} investigates optimality
conditions for the constrained QP~\eqref{EqnOrigQP}, and provides
deterministic conditions, of the mutual incoherence form, under which
a sparse solution, which is known to be within $\epsilon$ of the
observed values, can be recovered exactly.  Among a variety of other
results, both Tropp~\cite{Tropp06} and Donoho et al.~\cite{DonElTem06}
also provide sufficient conditions for the support of the optimal
solution to the constrained QP~\eqref{EqnOrigQP} to be contained
within the true support of $\betastar$.  Most directly related to the
current paper is recent work by both Meinshausen and
Buhlmann~\cite{Meinshausen06}, focusing on Gaussian noise, and
extensions by Zhao and Yu~\cite{Zhao06} to more general noise
distributions, on the use of the Lasso for model selection.  For the
case of Gaussian noise, both papers established that under mutual
incoherence conditions and appropriate choices of the regularization
parameter $\lambda_{\numobs}$, the Lasso can recover the sparsity
pattern with probability converging to one for particular regimes of
$\numobs$, $\mdim$ and $\spindex$, when $\arow_k$ drawn randomly from
random Gaussian ensembles. We discuss connections to our results at
more length in the the sequel.

\subsection{Our contributions}

Recall the linear observation model~\eqref{EqnLinearObs}.  For
compactness in notation, let us use $\Amat$ to denote the $\numobs
\times \mdim$ matrix formed with the vectors $\arow_k = (\arow_{k1},
\arow_{k2}, \ldots, \arow_{k \mdim}) \in \real^\mdim$ as rows, and the
vectors $\acol_j = (\arow_{1j}, \arow_{2j}, \ldots, \arow_{\numobs
j})^T \in \real^\numobs$ as columns, as follows:
\begin{eqnarray}
\label{EqnAmatDefn}
\Amat & \defn & \begin{bmatrix} \arow_1^T \\ \arow_2^T \\ \vdots \\
\arow_\numobs^T \end{bmatrix} \; = \; \begin{bmatrix} \acol_1 & \acol_2
& \cdots & \acol_\mdim \end{bmatrix}.
\end{eqnarray}
Consider the (random) set $\SolSet(\Amat, \betastar, \Wsca,
\lambda_\numobs)$ of optimal solutions to this constrained quadratic
program~\eqref{EqnOrigQP}.  By convexity and boundedness of the cost
function, the solution set is always non-empty.  For any vector
$\betasca \in \real^\mdim$, we define the sign function
\begin{eqnarray}
\sgn(\betasca_i) & \defn & \begin{cases} +1 & \mbox{if $\betasca_i >
    0$} \\ -1 & \mbox{if $\betasca_i < 0$} \\ 0 & \mbox{if $\betasca_i
    = 0$}.
		    \end{cases}
\end{eqnarray}
Of interest is the event that the Lasso~\eqref{EqnOrigQP} succeeds in
recovering the sparsity pattern of the unknown $\betastar$:
\begin{description}
\item[Property $\Prop(\Amat, \betastar, \Wsca, \lambda_\numobs)$:]
There exists an optimal solution $\betahat \in \SolSet(\Amat,
\betastar, \Wsca, \lambda_\numobs)$ with the property
\mbox{$\sgn(\betahat) = \sgn(\betastar)$.}
\end{description}
Our main result is that for a broad class of random Gaussian ensembles
based on covariance matrices satisfying mutual incoherence conditions,
there exist fixed constants $0 < \ThreshLow \leq 1$ and $1 \leq
\ThreshUp < +\infty$ such that for all $\threshbou > 0$, property
$\Prop(\Amat, \betastar, \Wsca, \lambda_\numobs)$ holds with high
probability (over the choice of noise vector $\Wsca$ and random matrix
$\Amat$) whenever
\begin{eqnarray}
\label{EqnThreshold}
\numobs & > & 2 (\ThreshUp + \threshbou) \, \spindex \, \log(\mdim -
\spindex) + \spindex + 1,
\end{eqnarray}
and \emph{conversely}, fails to hold with high probability whenever
\begin{eqnarray}
\label{EqnThresholdLOw}
\numobs & < & 2 (\ThreshLow - \threshbou) \, \spindex \, \log(\mdim -
\spindex) + \spindex + 1.
\end{eqnarray}
Moreover, for the special case of the uniform Gaussian ensemble (i.e.,
$\arow_k \sim N(0, I_{\mdim})$), we show that \mbox{$\ThreshLow =
\ThreshUp = 1$,} so that the threshold is sharp.  This threshold
result has a number of connections to previous work in the area that
focuses on special forms of scaling.  More specifically, as we discuss
in more detail in Section~\ref{SecComp}, in the special case of linear
scaling (i.e., $\numobs = \gamma \mdim$ for some $\gamma > 0$), this
theorem provides a noisy analog of results previously established for
basis pursuit in the noiseless case~\cite{Donoho04a,CandesTao05}.
Moreover, our result can also be adapted to an entirely different
scaling regime for $\numobs, \mdim$ and $\spindex$, as considered by a
separate body of recent work~\cite{Meinshausen06,Zhao06} on the
high-dimensional Lasso.

The remainder of this paper is organized as follows.  We begin in
Section~\ref{SecPrelim} with some necessary and sufficient conditions,
based on standard optimality conditions for convex programs, for
property $\Prop(\Amat, \betastar, \Wsca, \lambda_\numobs)$ to hold.
We then prove a consistency result for the case of deterministic
design matrices $\Amat$.  Section~\ref{SecMain} is devoted to the
statement and proof of our main result on the asymptotic behavior of
the lasso for random Gaussian ensembles.  We illustrate this result
via simulation in Section~\ref{SecExperimental}, and conclude with a
discussion in Section~\ref{SecDiscussion}.

\section{Some preliminary analysis}
\label{SecPrelim}

In this section, we provide necessary and sufficient conditions for
property $\Prop(\Amat, \betastar, \Wsca, \lambda_\numobs)$ to hold.
Based on these conditions, we then define collections of random
variables that play a central role in our analysis.  In particular,
the study of $\Prop(\Amat, \betastar, \Wsca, \lambda_\numobs)$ is
reduced to the study of the extreme order statistics of these random
variables.  We then state and prove a result about the behavior of the
Lasso for the case of a deterministic design matrix $\Amat$.

\subsection{Necessary and sufficient conditions}

We begin with a simple set of necessary and sufficient conditions for
property $\Prop(\Amat, \betastar, \Esca, \lambda_\numobs)$ to hold.
We note that this result is not essentially new (e.g.,
see~\cite{Fuchs04,Fuchs05,Meinshausen06,Tropp06,Zhao06} for variants),
and follows in a straightforward manner from optimality conditions for
convex programs~\cite{Hiriart1}; see Appendix~\ref{AppNecSuff} for
further details.  We define $\Sset \defn \{i \in \{1, \ldots, \mdim \}
\; \mid \; \betastar_i \neq 0 \}$ to be the support of $\betastar$,
and let $\Sbar$ be its complement.  For any subset $T \subseteq \{1,
2, \ldots, \mdim\}$, let $\Amatt{T}$ be the $\numobs \times |T|$
matrix with the vectors $\{\acol_i, i \in T\}$ as columns.
\blems
\label{LemNecSuffCond}
Assume that the matrix $\Amatt{\Sset}^T \Amatt{\Sset}$ is invertible.
Then, for any given $\lambda > 0$ and noise vector $\esca \in
\real^\numobs$, property $\Prop(\Amat, \betastar, \esca,
\lambda_\numobs)$ holds if and only if
\begin{subequations}
\begin{eqnarray}
\label{EqnPropA}
\left | \Amatt{\Sbar}^T \Amatt{\Sset} \left(\Amatt{\Sset}^T \Amatt{\Sset}
\right)^{-1} \left[ \frac{1}{n} \Amatt{\Sset}^T \esca - \lambda
\sgn(\betastar_\Sset) \right] - \frac{1}{\numobs} \Amatt{\Sbar}^T \esca
\right| & \leq & \lambda, \quad \mbox{and} \\
\label{EqnPropB}
\left | \betastar_\Sset + \left(\frac{1}{\numobs} \Amatt{\Sset}^T \Amatt{\Sset}
\right)^{-1} \left[\frac{1}{\numobs} \Amatt{\Sset}^T \esca - \lambda
\sgn(\betastar_\Sset) \right] \right| & > & 0,
\end{eqnarray}
\end{subequations}
where both of these vector inequalities should be taken elementwise.
\elems

For shorthand, define $\spec \defn \sgn(\betastar_\Sset)$, and denote
by $e_i \in \real^s$ the vector with $1$ in the $i^{th}$ position, and
zeroes elsewhere.  Motivated by Lemma~\ref{LemNecSuffCond}, much of
our analysis is based on the collections of random variables, defined
each index $i \in \Sset$ and $j \in \Sbar$ as follows:
\begin{subequations}
\begin{eqnarray}
\label{EqnDefnUvar}
\Uvar_i & \defn & e_i^T \left(\frac{1}{\numobs} \Amatt{\Sset}^T
\Amatt{\Sset} \right)^{-1} \left[\frac{1}{\numobs} \Amatt{\Sset}^T
\Wsca - \lambda_{\numobs} \spec \right] \\
\label{EqnDefnVvar}
\Vvar_j & \defn & \acol_j^T \Biggr \{ \Amatt{\Sset}
\left(\Amatt{\Sset}^T \Amatt{\Sset} \right)^{-1} \lambda_{\numobs}
\spec - \left[\Amatt{\Sset}\left(\Amatt{\Sset}^T \Amatt{\Sset}
\right)^{-1} \Amatt{\Sset}^T - I_{n \times n} \right]
\frac{\Esca}{\numobs} \Biggr \}.
\end{eqnarray}
\end{subequations}
Recall that $\spindex = |\Sset|$ and $\bigN = |\Sbar| = \mdim -
\spindex$.  From Lemma~\ref{LemNecSuffCond}, the behavior of
$\Prop(\Amat, \betastar, \Wsca, \lambda_\numobs)$ is determined by the
behavior of $\max_{j \in \Sbar} |\Vvar_j|$ and $\max_{i \in \Sset}
|\Uvar_i|$.  In particular, condition~\eqref{EqnPropA} holds if and
only if the event
\begin{eqnarray}
\label{EqnMaxV}
\evA & \defn & \left \{ \max_{j \in \Sbar} | \Vvar_j| \, \leq \,
\lambda_{\numobs} \right \}
\end{eqnarray}
holds.  On the other hand, if we define $\rho_n \defn \min_{i
\in \Sset} |\betastar_i|$, then the event
\begin{eqnarray}
\label{EqnMaxU}
\evB & \defn & \left \{ \max_{i \in \Sset} |\Uvar_i | \, \leq \,
\rho_n \right \}
\end{eqnarray}
is sufficient to guarantee that condition~\eqref{EqnPropB} holds.
Consequently, our proofs are based on analyzing the asymptotic
probability of these two events.

\subsection{Recovery of sparsity: deterministic design}

We now show how Lemma~\ref{LemNecSuffCond} can be used to analyze the
behavior of the Lasso for the special case of a deterministic
(non-random) design matrix $\Amat$.  To gain intuition for the
conditions in the theorem statement, it is helpful to consider the
\emph{zero-noise condition} $\esca = 0$, in which each observation
$\Ysca_k = \arow_k^T \betastar$ is uncorrupted.  In this case, the
conditions of Lemma~\ref{LemNecSuffCond} reduce to
\begin{subequations}
\begin{eqnarray}
\label{EqnPropASimp}
\left | \Amatt{\Sbar}^T \Amatt{\Sset} \left(\Amatt{\Sset}^T
\Amatt{\Sset}\right)^{-1} \sgn(\betastar_\Sset) \right | & \leq & 1 \\
\label{EqnPropBSimp}
\left | \betastar_\Sset - \lambda \left(\frac{1}{n} \Amatt{\Sset}^T
\Amatt{\Sset} \right)^{-1} \sgn(\betastar_\Sset) \right| & > & 0.
\end{eqnarray}
\end{subequations}
Of course, if the conditions of Lemma~\ref{LemNecSuffCond} fail to
hold in the zero-noise setting, then there is little hope of
succeeding in the presence of noise.

The zero-noise conditions motivate imposing the following set of
conditions on the design matrix:
\begin{subequations}
\label{EqnKeyCond}
\begin{eqnarray}
\label{EqnKeyCondA}
\left \| \Amatt{\Sbar}^T \Amatt{\Sset} \left(\Amatt{\Sset}^T
\Amatt{\Sset}\right)^{-1} \right \|_\infty & \leq & (1-\matbound) \quad
\mbox{for some $\matbound \in (0,1]$, and} \\
\label{EqnKeyCondB}
\myeigmin(\frac{1}{\numobs} \Amatt{\Sset}^T \Amatt{\Sset}) \geq \Cmin >
0,
\end{eqnarray}
\end{subequations}
where $\myeigmin$ denotes the minimal eigenvalue.  Under these
conditions, we have the following:
\bprops
\label{PropDetDesign}
Suppose that we observe $\Ysca = \Amat \betastar + \Esca$, where each
column $\acol_j$ of $\Amat$ is normalized to $\ell_2$-norm $\numobs$,
and $\Esca \sim N(0, \sigw I)$.  Assume $\betastar$ and $\Amat$
satisfy conditions~\eqref{EqnKeyCond}, and define \mbox{$\rho_\numobs
\defn \min_{i \in \Sset} |\betastar_i|$.}  If $\lambda_{\numobs}
\rightarrow 0$ is chosen such that
\begin{equation}
\mbox{(a)} \quad \frac{\numobs \lambda_{\numobs}^2}{\log
(\mdim-\spindex)} \rightarrow +\infty, \qquad \mbox{and} \quad
\mbox{(b)} \quad \frac{1}{\rho_n} \; \left \{ \sqrt{\frac{\log
\spindex}{\numobs }} + \lambda_\numobs \, \| (\frac{1}{\numobs}
\Amatt{\Sset}^T \Amatt{\Sset})^{-1} \|_\infty \right \} \rightarrow 0,
\end{equation}
then $\Prob(\Prop(\Amat, \betastar, \Esca, \lambda_\numobs)
\rightarrow 1$ as $\numobs \rightarrow +\infty$.
\eprops

Before proving the proposition, we pause to make a number of comments.
First, conditions of the form~\eqref{EqnKeyCondA} have been considered
in previous work on the
lasso~\cite{Fuchs04,Fuchs05,Meinshausen06,Tropp06,Zhao06}.  In
particular, various authors~\cite{Tropp06,Meinshausen06,Zhao06}
provide examples and results on matrix families that satisfy this type
of condition.  Moreover, previous work~\cite{Meinshausen06,Zhao06}
provides asymptotic results for particular scalings of $\mdim$,
$\spindex$ and $\numobs$ for random design matrices, as we discuss in
more detail in Section~\ref{SecMain}.  To the best of our knowledge,
Proposition~\ref{PropDetDesign} is the first result to provide
sufficient conditions for exact recovery in deterministic designs with
general scaling of $\mdim$, $\spindex$ and $\numobs$.

Second, it is worthwhile to consider Proposition~\ref{PropDetDesign}
in the classical setting (i.e., in which the number of samples
$\numobs \rightarrow +\infty$ with $\mdim$ and $\spindex$ remaining
fixed).  In this setting, the quantity $\rho_\numobs = \min_{i \in
\Sset} |\betastar_i|$ does not depend on $\numobs$.  Hence, in
addition to the condition~\eqref{EqnKeyCond}, the requirements reduce
to $\lambda_{\numobs} \rightarrow 0$ and $ n \lambda_{\numobs}^2
\rightarrow +\infty$.  Note that $\lambda_{\numobs} = \frac{ \log
\numobs }{\sqrt{\numobs}}$ is one suitable choice.  This classical
case is also covered by previous work~\cite{Knight00,Meinshausen06,Zhao06}.

Last, consider the more general setting where all three parameters
$(\numobs, \mdim, \spindex)$ grow to infinity, and suppose for
simplicity that $\rho_\numobs$ stays bounded away from $0$.  The
conditions $\lambda^2_\numobs \rightarrow 0$ and $\lambda_\numobs^2 \;
\frac{\numobs}{\log (\mdim - \spindex)} \rightarrow +\infty$ imply
that the number of observations $\numobs$ must grow at a rate faster
than $\log (\mdim - \spindex)$.  In the following section, in which we
consider the more general case of random Gaussian ensembles, we will
see that for ensembles satisfying mutual incoherence conditions, we in
fact require that $\frac{\numobs}{\log(\mdim-\spindex)} =
\Theta(\spindex) \rightarrow +\infty$.

\vtiny

\subsection{Proof of Proposition~\ref{PropDetDesign}}

Recall the events $\evA$ and $\evB$ defined in
equations~\eqref{EqnMaxV} and~\eqref{EqnMaxU} respectively.  To
establish the claim, we must show that that $\Prob[\evA^c \; \mbox{or}
\; \evB^c] \rightarrow 0$, where $\evA^c$ and $\evB^c$ denote the
complements of these events. By union bound, it suffices to show both
$\Prob[\evA^c]$ and $\Prob[\evB^c]$ converge to zero, or equivalently
that $\Prob[\evA]$ and $\Prob[\evB]$ both converge to one.

\noindent \paragraph{Analysis of $\evA$:} We begin by establishing
that $\Prob[\evA] \rightarrow 1$. Throughout the proof, we use the
shorthand $\svec \defn \sgn(\betastar)$ and $\bigN \defn \mdim - \spindex = |\Sbar|$.

Recalling the definition~\eqref{EqnDefnVvar} of the random variables
$\Vvar_j$, note that $\evA$ holds holds if and only $\frac{\min_{j \in
\Sbar} \Vvar_j }{\lambda_{\numobs}} \geq -1$ and $\frac{\max_{j \in
\Sbar} \Vvar_j}{\lambda_{\numobs}} \leq 1$.  Moreover, we note that
each $\Vvar_j$ is Gaussian with mean
\begin{eqnarray*}
\mu_j \; = \; \Exs[\Vvar_j] & = & \lambda_{\numobs} \acol_j^T
\Amatt{\Sset} \left(\Amatt{\Sset}^T \Amatt{\Sset} \right)^{-1} \svec.
\end{eqnarray*}
Using condition~\eqref{EqnKeyCondA}, we have $|\mu_j| \leq
(1-\matbound) \, \lambda_{\numobs}$ for all indices $j=1,\ldots,
\bigN$, from which we obtain that
\begin{equation*}
\frac{\max_{j \in \Sbar} \Vvar_j}{\lambda_{\numobs}} \: \leq \:
(1-\matbound) + \frac{1}{\lambda_{\numobs}} \max_{j} \Vtil_j, \qquad
\mbox{and} \qquad
\frac{\min_{j \in \Sbar} \Vvar_j}{\lambda_{\numobs}} \: \geq |;
-(1-\matbound) + \frac{1}{\lambda_{\numobs}} \min_{j} \Vtil_j,
\end{equation*}
where $\Vtil_j \defn \acol_j^T \left[I_{\numobs \times \numobs} -
\Amatt{\Sset} \left(\Amatt{\Sset}^T \Amatt{\Sset} \right)^{-1}
\Amatt{\Sset}^T \right] \Esca$ are zero-mean (correlated) Gaussian
variables.  Hence, in order to establish condition~\eqref{EqnPropA} of
Lemma~\ref{LemNecSuffCond}, we need to show that
\begin{equation}
\Prob \left[\frac{1}{\lambda_{\numobs}} \min_{j \in \Sbar} \Vtil_j <
  -\matbound, \quad \mbox{or} \quad \frac{1}{\lambda_{\numobs}}
  \max_{j \in \Sbar} \Vtil_j > \matbound \right] \rightarrow 0.
\end{equation}
In fact, using Lemma~\ref{LemSimpleInequal} (see
Appendix~\ref{AppAux}), it is sufficient to show that
\mbox{$\Prob[\frac{\max_{j \in \Sbar} |\Vtil_j|} {\lambda_{\numobs}} >
\matbound] \rightarrow 0$}.  By applying Markov's inequality and
Gaussian comparison results~\cite{LedTal91} (see
Lemma~\ref{LemGeneric} in Appendix~\ref{AppGaussComp}), we obtain
\begin{equation*}
\Prob \left[\frac{\max_{j \in \Sbar} |\Vtil_j|} {\lambda_{\numobs}} >
  \matbound\right] \; \leq \; \frac{\Exs[\max_{j \in \Sbar}
  |\Vtil_j|]}{\lambda_{\numobs}} \; \leq \; \frac{3 \sqrt{\log N}
  }{\lambda_{\numobs}} \max_j \sqrt{\Exs [\Vtil_j^2]}.
\end{equation*}
Straightforward computation yields that
\[
\Exs[\Vtil_j^2] = \frac{\sigw}{\numobs^2} \; \acol_j^T \left[
    I_{\numobs \times \numobs} - \Amatt{\Sset} \left(\Amatt{\Sset}^T
    \Amatt{\Sset} \right)^{-1} \Amatt{\Sset}^T \right] \acol_j \; \leq
    \frac{\sigw}{\numobs^2} \|\acol_j\|^2 \; = \;
    \frac{\sigw}{\numobs},
\]
since the matrix $I_{\numobs \times \numobs} - \Amatt{\Sset}
\left(\Amatt{\Sset}^T \Amatt{\Sset} \right)^{-1} \Amatt{\Sset}^T$ has
maximum eigenvalue equal to one, and $\|\acol_j\|_2^2 = n$ by
construction.  Consequently, condition (a) in the theorem
statement---namely, that $\frac{\log \bigN}{\numobs \lambda^2_\numobs}
\rightarrow 0$ is sufficient to ensure that
$\Exs[\Vtilmax]/\lambda_{\numobs} \rightarrow 0$.  Thus, we have
established $\Prob(\evA) \rightarrow 1$ (i.e., that
condition~\eqref{EqnPropA} holds w.p. one as $\numobs \rightarrow
+\infty$).

\vtiny

\noindent \paragraph{Analysis of $\evB$:} We now show that
$\Prob(\evB) \rightarrow 1$.  Beginning with the triangle inequality,
we upper bound $\max_{i} |\Uvar_i| \defn \| (\frac{1}{\numobs}
\Amatt{\Sset}^T \Amatt{\Sset} )^{-1} [\frac{1}{\numobs}
\Amatt{\Sset}^T \Wsca -\lambda_{\numobs} \sgn(\betastar_\Sset) ]
\|_\infty$ as
\begin{eqnarray*}
\max_{i} |\Uvar_i| & \leq & \left \| (\frac{1}{\numobs}
\Amatt{\Sset}^T \Amatt{\Sset} )^{-1} \frac{1}{\numobs} \Amatt{\Sset}^T
\Wsca \right \|_\infty + \left \| (\frac{1}{\numobs} \Amatt{\Sset}^T
\Amatt{\Sset})^{-1} \right \|_\infty \; \lambda_{\numobs} 
\end{eqnarray*}
Let $e_i$ denote the unit vector with one in position $i$ and zeroes
elsewhere.  Now define, for each index $i \in \Sset$, the Gaussian
random variable $\Zvar_i \defn e_i^T (\frac{1}{\numobs}
\Amatt{\Sset}^T \Amatt{\Sset} )^{-1} \frac{1}{\numobs} \Amatt{\Sset}^T
\Wsca$.  Each such $\Zvar_i$ is a zero-mean Gaussian with variance
given by
\begin{eqnarray*}
\var(\Zvar_i) & = & \frac{\sigw}{\numobs} e_i^T (\frac{1}{\numobs}
\Amatt{\Sset}^T \Amatt{\Sset})^{-1} e_i \; \leq \; \frac{\sigw}{\Cmin
\numobs}
\end{eqnarray*}
Hence, by a standard Gaussian comparison theorem~\cite{LedTal91} (in
particular, see Lemma~\ref{LemGeneric} in
Appendix~\ref{AppGaussComp}), we have
\begin{eqnarray*}
\Exs [\max_{1 \leq i \leq \spindex} |\Zvar_i|] & = & \Exs \left[\left
\| (\frac{1}{\numobs} \Amatt{\Sset}^T \Amatt{\Sset} )^{-1}
\frac{1}{\numobs} \Amatt{\Sset}^T \Wsca \right \|_\infty \right] \\
& \leq & 3 \sqrt{\frac{\sigw \log \spindex}{\numobs \Cmin}}.
\end{eqnarray*}
Thus, recalling the defining $\rho_n \defn \min_{i \in \Sset}
|\betastar_i|$, we apply Markov's inequality to conclude that
\begin{eqnarray*}
1 - \Prob \left[\left | \betastar_\Sset + \left(\frac{1}{\numobs}
\Amatt{\Sset}^T \Amatt{\Sset} \right)^{-1} \left[\frac{1}{\numobs}
\Amatt{\Sset}^T \esca - \lambda \sgn(\betastar_\Sset) \right] \right|
> 0 \right] & \leq & \Prob \left[ \frac{1}{\rho_n} \max_{1 \leq i \leq
\spindex} |\Uvar_i| > 1 \right] \\
& \leq & \Prob \left[ \frac{1}{\rho_n} \left \{ \max_{1 \leq i \leq
\spindex} |\Zvar_i| + \lambda_{\numobs} \| (\frac{1}{\numobs}
\Amatt{\Sset}^T \Amatt{\Sset})^{-1} \|_\infty \right \} > 1 \right]
\\
& \leq & \frac{1}{\rho_n} \left \{ \Exs \left[\max_{1 \leq i \leq
    \spindex} |\Zvar_i|\right] +  \lambda_{\numobs} \|
    (\frac{1}{\numobs} \Amatt{\Sset}^T \Amatt{\Sset})^{-1} \|_\infty \right \}
    \\
& \leq & \frac{1}{\rho_n} \; \left \{ 3 \sqrt{\frac{\sigw \log
\spindex}{\numobs \Cmin}} + \lambda_{\numobs} \| (\frac{1}{\numobs}
\Amatt{\Sset}^T \Amatt{\Sset})^{-1} \|_\infty \right \},
\end{eqnarray*}
which converges to zero as $\numobs \rightarrow +\infty$, using
condition (b) in the theorem statement.
\hfill \qed

\section{Recovery of sparsity: random Gaussian ensembles}
\label{SecMain}

We now turn to the analysis of random design matrices $\Amat$, in
which each row $\arow_k$ is chosen as an i.i.d. Gaussian random vector
with covariance matrix $\CovMat$.  In particular, we prove the
existence of thresholds that provide a sharp description of the
failure/success of the Lasso as a function of $(\numobs, \mdim,
\spindex)$.  We begin by setting up and providing a precise statement
of the main result, and then discussing its connections to previous
work.  In the later part of this section, we provide the proof.

\subsection{Statement of main result}

Consider a covariance matrix $\CovMat$ with unit diagonal, and with
its minimum and maximum eigenvalues (denoted $\myeigmin$ and
$\myeigmax$ respectively) bounded as
\begin{equation}
\label{EqnEigCond}
\myeigmin(\CovMat_{\Sset \Sset}) \geq \Cmin, \qquad \mbox{and} \qquad
\myeigmax(\CovMat) \leq \Cmax
\end{equation}
for constants $\Cmin > 0$ and $\Cmax < +\infty$.  Given a vector
$\betastar \in \real^\mdim$, define its support $\Sset = \{ i \in \{1,
\ldots, \mdim\} \; | \; \betastar_i \neq 0 \}$, as well as the
complement $\Sbar$ of its support.  Suppose that $\CovMat$ and $\Sset$
satisfy the conditions $\| (\CovMat_{\Sset \Sset})^{-1} \|_\infty \leq
\Dconmax$ for some $\Dconmax < +\infty$, and
\begin{eqnarray}
\label{EqnCovCond}
\| \CovMat_{\Sbar \Sset} (\CovMat_{\Sset\Sset})^{-1} \|_\infty & \leq
& (1-\matbound)
\end{eqnarray}
for some $\matbound \in (0,1]$.  Under these conditions, we consider
the observation model
\begin{eqnarray}
\label{EqnObsModel}
\Ysca_k & = & \arow_k^T \betastar + \Esca_k, \qquad k = 1, \ldots,
\numobs,
\end{eqnarray}
where $\arow_k \sim N(0, \CovMat)$ and $\Esca_k \sim N(0, \sigma^2)$
are independent Gaussian variables for $k= 1, \ldots, \numobs$.
Furthermore, we define $\rho_\numobs \defn \min_{i \in \Sset}
|\betastar_i|$, and the sparsity index $\spindex = |\Sset|$.

\btheos
\label{ThmGenNoiseGauss}
Consider a sequence of covariance matrices $\{\CovMat[\mdim] \}$ and
solution vectors $\{\betastar[\mdim]\}$ satisfying
conditions~\eqref{EqnEigCond} and~\eqref{EqnCovCond}.  Under the
observation model~\eqref{EqnObsModel}, consider a sequence $(\numobs,
\mdim(\numobs), \spindex(\numobs))$ such that $\spindex$, $(\numobs -
\spindex)$ and $(\mdim-\spindex)$ tend to infinity. Define the
thresholds
\begin{equation}
\ThreshLow \, \defn \, \frac{(\sqrt{\Cmax} - \sqrt{\Cmax -
\frac{1}{\Cmax}})^2}{\Cmax \, (2-\matbound)^2} \leq 1, \qquad
\mbox{and} \qquad \ThreshUp \, \defn \, \frac{\Cmax}{\matbound^2
\Cmin} \geq 1.
\end{equation}
Then for any constant $\threshbou > 0$, we have the following
\begin{enumerate}
\item[(a)] If $\numobs < 2 (\ThreshLow - \threshbou) \, \spindex \,
\log(\mdim-\spindex) + \spindex + 1$, then $\Prob[\Prop(\Amat,
\betastar, \Esca, \lambda_\numobs)] \rightarrow 0$ for any
non-increasing sequence $\lambda_{\numobs} > 0$.
\item[(b)] Conversely, if $\numobs > 2 (\ThreshUp + \threshbou) \,
\spindex \, \log(\mdim-\spindex) + \spindex$, and $\lambda_{\numobs}
\rightarrow 0$ is chosen such that
\begin{equation}
\label{EqnLamGaussCondition}
\frac{\numobs \lambda^2_\numobs}{ \log (\mdim-\spindex)} \rightarrow
+\infty, \qquad \mbox{and} \qquad
\frac{1}{\rho_\numobs} \Big[ \lambda_{\numobs} + \sqrt{ \frac{\log
  \spindex}{\numobs}} \Big] \; \rightarrow \; 0,
\end{equation}
then $\Prob[\Prop(\Amat, \betastar, \Esca, \lambda_\numobs)]
\rightarrow 1$.
\end{enumerate}

\etheos
\noindent {\bf{Remark:}} Suppose for simplicity that $\rho_\numobs$
remains bounded away from $0$.  In this case, the requirements on
$\lambda_\numobs$ reduce to $\lambda_\numobs \rightarrow 0$, and
$\lambda^2_\numobs \numobs/\log(\mdim - \spindex) \rightarrow
+\infty$.  One suitable choice is $\lambda^2_\numobs = \frac{
\log(\spindex) \; \log (\mdim - \spindex) }{\numobs}$, with which we
have
\begin{eqnarray*}
\lambda^2_\numobs & = & \left (\frac{ \spindex \, \log
(\mdim-\spindex)} {\numobs} \right) \; \frac{\log(\spindex)}{\spindex}
\; = \; O \left (\frac{\log \spindex}{\spindex} \right) \; \rightarrow
\; 0,
\end{eqnarray*}
and
\begin{eqnarray*}
\frac{\numobs \lambda_\numobs^2}{\log (\mdim - \spindex)} & = & \log
(\spindex) \rightarrow +\infty.
\end{eqnarray*}
Without a bound on $\rho_\numobs$, the second condition in
equation~\eqref{EqnLamGaussCondition} constrains the rate of decrease
of the minimum \mbox{$\rho_\numobs = \min_{i \in \Sset}
|\betastar_i|$.}

\subsection{Some consequences}
\label{SecComp}

To develop intuition for this result, we begin by stating certain
special cases as corollaries, and discussing connections to previous
work.

\subsubsection{Uniform Gaussian ensembles}

First, we consider the special case of the uniform Gaussian ensemble,
in which $\CovMat = I_{\mdim \times \mdim}$.  Previous work by
Donoho~\cite{Donoho04a} as well as Candes and Tao~\cite{CandesTao05}
has focused on the uniform Gaussian ensemble in the the noiseless
($\sigma^2 = 0$) and underdetermined setting ($\numobs = \gamma \mdim$
for some $\gamma \in (0,1)$).  Analyzing the asymptotic behavior of
the linear program~\eqref{EqnLPRelax} for recovering $\betastar$, the
basic result is that there exists some $\alpha > 0$ such that all
sparsity patterns with $\spindex \leq \alpha \mdim$ can be recovered
with high probability.

Applying Theorem~\ref{ThmGenNoiseGauss} to the noisy version of this
problem, the uniform Gaussian ensemble means that we can choose
$\matbound = 1$, and $\Cmin = \Cmax = 1$, so that the threshold
constants reduce
\begin{equation*}
\ThreshLow \, = \, \frac{(\sqrt{\Cmax} - \sqrt{\Cmax -
\frac{1}{\Cmax}})^2}{\Cmax \, (2-\matbound)^2} \, = \, 1 \qquad
\mbox{and} \qquad
\ThreshUp \, = \, \frac{\Cmax}{\matbound^2 \Cmin} \, = \, 1.
\end{equation*}
Consequently, Theorem~\ref{ThmGenNoiseGauss} provides a sharp
threshold for the behavior of the Lasso, in that failure/success is
entirely determined by whether or not $\numobs > 2 \spindex \, \log
(\mdim-\spindex) + \spindex + 1$.  Thus, if we consider the particular
linear scaling analyzed in previous work on the noiseless
case~\cite{Donoho04a,CandesTao05}, we have:
\bcors[Linearly underdetermined setting] 
\label{Cor1}
Suppose that $\numobs = \gamma \mdim$ for some $\gamma \in (0,1)$.  Then
\begin{enumerate}
\item[(a)] If $\spindex = \alpha \mdim$ for any $\alpha \in (0, 1)$,
then $\Prob\left[ \Prop(\Amat, \betastar, \Esca,
\lambda_\numobs)\right] \rightarrow 0$ for any positive sequence
$\lambda_{\numobs} > 0$.
\item[(b)] On the other hand, if $\spindex = O( \frac{\mdim}{\log
\mdim})$, then $\Prob\left[ \Prop(\Amat, \betastar, \Esca,
\lambda_\numobs)\right] \rightarrow 1$ for any sequence
$\{\lambda_{\numobs}\}$ satisfying the conditions of
Theorem~\ref{ThmGenNoiseGauss}(a).
\end{enumerate}
\ecors
\noindent Conversely, suppose that the size $\spindex$ of the support
of $\betastar$ scales linearly with the number of parameters $\mdim$.
The following result describes the amount of data required for the
$\ell_1$-constrained QP to recover the sparsity pattern in the noisy
setting ($\sigma^2 > 0$):
\bcors[Linear fraction support]
Suppose that $\spindex = \alpha \mdim$ for some $\alpha \in (0,1)$.
Then we require $\numobs > 2 \alpha \mdim \log [(1-\alpha) \, \mdim] +
\alpha \mdim $ in order to obtain exact recovery with probability
converging to one for large problems.
\ecors
\noindent These two corollaries establish that there is a significant
difference between recovery using basis pursuit~\eqref{EqnLPRelax} in
the noiseless setting versus recovery using the
Lasso~\eqref{EqnOrigQP} in the noisy setting.  When the amount of data
$\numobs$ scales only linearly with ambient dimension $\mdim$, then
the presence of noise means that the recoverable support size drops
from a linear fraction (i.e., $\spindex = \alpha \mdim$ as in the
work~\cite{Donoho04a,CandesTao05}) to a sublinear fraction (i.e.,
$\spindex = O(\frac{\log \mdim}{\mdim})$, as in Corollary~\ref{Cor1}).

\subsubsection{Non-uniform Gaussian ensembles}

We now consider more general (non-uniform) Gaussian ensembles that
satisfy conditions~\eqref{EqnEigCond} and~\eqref{EqnCovCond}.  As
mentioned earlier, previous papers by both Meinshausen and
Buhlmann~\cite{Meinshausen06} as well as Zhao and Yu~\cite{Zhao06}
treat model selection with the high-dimensional Lasso.  For suitable
covariance matrices (e.g., satisfying conditions~\eqref{EqnEigCond}
and~\eqref{EqnCovCond}), both sets of authors proved that the sparsity
pattern can be recovered exactly under scaling conditions of the form
\begin{equation}
\label{EqnMeinshausen}
\spindex = O(\numobs^{c_1}), \qquad \mbox{and} \quad\mdim =
O(e^{\numobs^{c_2}}), \qquad \mbox{where} \quad c_1 + c_2 < 1.
\end{equation}
Applying Theorem~\ref{ThmGenNoiseGauss} in this scenario, we have the
following:
\bcors Under the scaling~\eqref{EqnMeinshausen}, the Lasso will
recover the sparsity pattern with probability converging to one.
\ecors
\spro
Substituting the conditions~\eqref{EqnMeinshausen} into the threshold
condition~\eqref{EqnThreshold}, we obtain that the RHS takes the form
\begin{eqnarray*}
2 \spindex \log (\mdim - \spindex) + \spindex + 1 & = &
O(\numobs^{c_1}) \log \left[ O(e^{\numobs^{c_2}}) - O(\numobs^{c_1})
\right] + O(\numobs^{c_1}) \\
& = & O(\numobs^{c_1 +c_2}) \; \ll \; \numobs,
\end{eqnarray*}
since $c_1 + c_2 < 1$ by assumption. Thus, we see that under these
conditions, our threshold condition~\eqref{EqnThreshold} is satisfied
\emph{a fortiori}.  
\fpro
\noindent In fact, under this stronger scaling~\eqref{EqnMeinshausen},
both papers~\cite{Meinshausen06,Zhao06} proved that the probability of
exact recovery converges to one at a rate exponential in some
polynomial function of $\numobs$.  Interestingly, our results show
that the Lasso can recover the sparsity pattern for a much broader
range of $(\numobs, \mdim, \spindex)$ scaling.

\subsection{Proof of Theorem~\ref{ThmGenNoiseGauss}(b)}

We now turn to the proof of part (b) of our main result.  As with the
proof of Proposition~\ref{PropDetDesign}, the proof is based on
analyzing the collections of random variables $\{ \Vvar_j \; \mid \; j
\in \Sbar \}$ and $\{ \Uvar_i \; \mid \; i \in \Sset \}$, as defined
in equations~\eqref{EqnDefnUvar} and~\eqref{EqnDefnVvar} respectively.
We begin with some preliminary results that serve to set up the
argument.

\subsubsection{Some preliminary results}

We first note that for $\spindex < \numobs$, the random Gaussian
matrix $\Amatt{\Sset}$ will have rank $\spindex$ with probability one,
whence the matrix $\Amatt{\Sset}^T \Amatt{\Sset}$ is invertible with
probability one.  Accordingly, the necessary and sufficient conditions
of Lemma~\ref{LemNecSuffCond} are applicable.  Our first lemma, proved
in Appendix~\ref{AppGaussCond}, concerns the behavior of the random
vector $\Vvar = (\Vvar_1, \ldots, \Vvar_\bigN)$, when conditioned on
$\Amatt{\Sset}$ and $\Esca$.  Recalling the shorthand notation $\spec
\defn \sgn(\betastar)$, we summarize in the following
\blems
\label{LemGaussCond}
Conditioned on $\Amatt{\Sset}$ and $\Esca$, the random vector $(\Vvar
\, \mid \, \Esca, \Amatt{\Sset})$ is Gaussian.  Its mean vector is
upper bounded as
\begin{eqnarray}
\label{EqnMeanBound}
\left | \Exs[\Vvar \, \mid \, \Esca, \Amatt{\Sset}] \right | & \leq &
\lambda_{\numobs} (1-\matbound) \, \ones.
\end{eqnarray}
Moreover, its conditional covariance takes the form
\begin{eqnarray}
\label{EqnVCov}
\cov[\Vvar\, \mid \, \Esca, \Amatt{\Sset}] & = & \myVar
\CovMat_{(\Sbar \, | \, \Sset)} \; = \; \myVar \, \big[\CovMat_{\Sbar
    \Sbar} - \CovMat_{\Sbar \Sset} (\CovMat_{\Sset\Sset})^{-1}
  \CovMat_{\Sset \Sbar} \big],
\end{eqnarray}
where
\begin{eqnarray}
\label{EqnDefnMyVar}
\myVar & \defn & \lambda_{\numobs}^2 \spec^T (\Amatt{\Sset}^T \Amatt{\Sset})^{-1}
\spec + \frac{1}{n^2} \Esca^T \left[I_{n \times n} - \Amatt{\Sset}
\left(\Amatt{\Sset}^T \Amatt{\Sset} \right)^{-1} \Amatt{\Sset}^T \right] \Esca
\end{eqnarray}
is a random scaling factor.
\elems

The following lemma, proved in Appendix~\ref{AppExtStats}, captures
the behavior of the random scaling factor $\myVar$ defined in
equation~\eqref{EqnDefnMyVar}:
\blems
\label{LemExtStats}
The random variable $\myVar$ has mean
\begin{eqnarray}
\Exs[\myVar] & = & \frac{\lambda_{\numobs}^2}{\numobs-\spindex-1} \,
\spec^T (\CovMat_{\Sset\Sset})^{-1} \spec + \frac{\sigw \,
(\numobs-\spindex)}{\numobs^2}.
\end{eqnarray}
Moreover, it is sharply concentrated in that for any $\delta > 0$, we
have
\begin{eqnarray}
\Prob \left[ \big | \myVar - \Exs[\myVar] \big | \geq \delta
\Exs[\myVar] \right] & \rightarrow & 0 \qquad \qquad \mbox{as $\numobs
\rightarrow +\infty$.}
\end{eqnarray}

\elems

\subsubsection{Main argument} 

With these preliminary results in hand, we now turn to analysis of the
collections of random variables $\{\Uvar_i, i \in \Sset\}$ and
$\{\Vvar_j, j \in \Sbar\}$.

\paragraph{Analysis of $\evA$:}
We begin by analyzing the behavior of $\max_{j \in \Sbar} |V_j|$.
First, for a fixed but arbitrary $\delta > 0$, define the event
$\Tail(\delta) \defn \{ |\myVar - \Exs[\myVar]| \geq \delta
\Exs[\myVar] \}$.  By conditioning on $\Tail(\delta)$ and its
complement $[\Tail(\delta)]^c$, we have the upper bound
\begin{eqnarray*}
\Prob[\max_{j \in \Sbar} |\Vvar_j| > \lambda_{\numobs}] & \leq & \Prob
\left[\max_{j \in \Sbar} |\Vvar_j| > \lambda_{\numobs} \, \mid \,
[\Tail(\delta)]^c \right] + \Prob[\Tail(\delta)].
\end{eqnarray*}
By the concentration statement in Lemma~\ref{LemExtStats}, we have
$\Prob[\Tail(\delta)] \rightarrow 0$, so that it suffices to analyze
the first term.  Set $\mu_j = \Exs[\Vvar_j | \Amatt{\Sset}]$, and let
$\Zvar$ be a zero-mean Gaussian vector with $\cov(Z) = \cov(\Vvar \, |
\, \Amatt{\Sset}, \Wsca)$.
\begin{eqnarray*}
\max_{j \in \Sbar} |\Vvar_j| & = & \max_{j \in \Sbar} |\mu_j +
\Zvar_j| \\
& \leq & \max_{j \in \Sbar} \left[ |\mu_j| + |\Zvar_j| \right] \\
& \leq & (1-\matbound) \lambda_{\numobs} + \max_{j \in \Sbar} |\Zvar_j|,
\end{eqnarray*}
where we have used the upper bound~\eqref{EqnMeanBound} on the
mean.  This inequality establishes the inclusion of events
\begin{eqnarray*}
\{\max_{j \in \Sbar}|\Zvar_j| \leq \matbound \lambda_{\numobs} \} &
\subseteq & \{\max_{j \in \Sbar} |\Vvar_j| \leq \lambda_{\numobs} \},
\end{eqnarray*}
thereby showing that it suffices to prove that $\Prob[\max_{j \in
\Sbar}|\Zvar_j| > \matbound \lambda_{\numobs} \, \mid \, [\Tail(\delta)]^c]
\rightarrow 0$.

Note that conditioned on $[\Tail(\delta)]^c$, the maximum value of
$\myVar$ is $\vstar \defn (1+\delta) \Exs[\myVar]$.  Since Gaussian
maxima increase with increasing variance, we have
\begin{eqnarray*}
\Prob \left[\max_{j \in \Sbar} |\Zvar_j| > \matbound \lambda_{\numobs}
\, \mid \, [\Tail(\delta)]^c \right] & \leq & \Prob \left[\max_{j \in
\Sbar} |\Zvarcond_j| > \matbound \lambda_{\numobs} \right],
\end{eqnarray*}
where $\Zvarcond$ is zero-mean Gaussian with covariance $\vstar \,
\CondMat$.

Using Lemma~\ref{LemSimpleInequal}, it suffices to show that $\Prob [
\max_{j \in \Sbar} \Zvarcond_j > \matbound \lambda_{\numobs}]$
converges to zero.  Accordingly, we complete this part of the proof
via the following two lemmas, both of which are proved in
Appendix~\ref{AppOneTwo}:
\blems
\label{LemOne}
Under the stated assumptions of the theorem, we have
$\frac{\vstar}{\lambda_{\numobs}^2} \rightarrow 0$ and
\begin{eqnarray*}
\lim_{n \rightarrow +\infty} \frac{1}{\lambda_{\numobs}} \Exs[\max_{j
\in \Sbar} \Zvarcond_j ] & \leq & \matbound.
\end{eqnarray*}
\elems
\blems
\label{LemTwo}
For any $\eta > 0$, we have
\begin{eqnarray}
\label{EqnLipCon}
\Prob\left[ \max_{j \in \Sbar} \Zvarcond_j > \eta + \Exs[\max_{j \in
\Sbar} \Zvarcond_j ] \right] & \leq & \exp \left (-\frac{\eta^2}{2 \vstar}
\right).
\end{eqnarray}
\elems
Lemma~\ref{LemOne} implies that for all $\delta > 0$, we have
$\Exs[\max_{j \in \Sbar} \Zvarcond_j ] \leq (1+\frac{\delta}{2})
\matbound \lambda_{\numobs}$ for all $n$ sufficiently large.
Therefore, setting $\eta = \frac{\delta}{2} \lambda_{\numobs}
\matbound$ in the bound~\eqref{EqnLipCon}, we have for fixed $\delta >
0$ and $n$ sufficiently large:
\begin{eqnarray*}
\Prob \left[ \max_{j \in \Sbar} \Zvarcond_j > (1+\delta)
\lambda_{\numobs} \matbound \right ] & \leq & \Prob \left[ \max_{j \in
\Sbar} \Zvarcond_j > \frac{\delta}{2} \lambda_{\numobs} \matbound +
\Exs[\max_{j \in \Sbar} \Zvarcond_j] \right] \\
& \leq & 2 \exp \left(- \frac{\delta^2 \lambda_{\numobs}^2
\matbound^2}{8 \vstar} \right).
\end{eqnarray*}
From Lemma~\ref{LemOne}, we have $\lambda_{\numobs}^2/\vstar
\rightarrow +\infty$, which implies that $\Prob [ \max_{j \in \Sbar}
\Zvarcond_j > (1+\delta) \lambda_{\numobs} \matbound ] \rightarrow 0$
for all $\delta > 0$.  By the arbitrariness of $\delta > 0$, we thus
have $\Prob[\max_{j \in \Sbar} \Zvarcond_j \leq \matbound
\lambda_{\numobs}] \rightarrow 1$, thereby establishing that
property~\eqref{EqnPropA} of Lemma~\ref{LemNecSuffCond} holds w.p. one
asymptotically.

\vtiny

\paragraph{Analysis of $\{\Uvar_i \}$:}
Next we prove that $\max_{i \in \Sset} |\Uvar_i| < \rho_\numobs \defn
\min_{i \in \Sset} |\betastar_i|$ with probability one as $\numobs
\rightarrow +\infty$.  Conditioned on $\Amatt{\Sset}$, the only random
component in $\Uvar_i$ is the noise vector $\Esca$.  A straightforward
calculation yields that this conditioned RV is Gaussian, with mean and
variance
\begin{eqnarray*}
\myumean_i \; \defn \; \Exs[\Uvar_i \; \mid \; \Amatt{\Sset}] & = &
-\lambda_{\numobs} e_i^T \left(\frac{1}{\numobs} \Amatt{\Sset}^T
\Amatt{\Sset} \right)^{-1} \spec, \\
\myuvar_i \; \defn \; \var[\Uvar_i \; \mid \; \Amatt{\Sset}] & = &
\frac{\sigw}{\numobs} e_i^T 
\left[ \frac{1}{\numobs} \Amatt{\Sset}^T
\Amatt{\Sset} \right]^{-1} e_i,
\end{eqnarray*}
respectively.  The following lemma, proved in
Appendix~\ref{AppUbehave}, is key to our proof:
\blems
\label{LemUbehave}
\noindent (a) The random variables $\myumean_i$ and $\myuvar_i$ have
means
\begin{equation}
\Exs[\myumean_i] = \frac{-\lambda_\numobs \;
\numobs}{\numobs-\spindex-1} e_i^T \; (\CovMat_{\Sset\Sset})^{-1} \,
\spec, \qquad \mbox{and} \qquad \Exs[\myuvar_i] \; = \;
\frac{\sigw}{\numobs-\spindex-1} \; e_i^T (\CovMat_{\Sset\Sset})^{-1}
e_i,
\end{equation}
respectively, which are bounded as
\begin{equation}
\label{EqnUbehaveBounds} 
|\Exs[\myumean_i]| \: \leq \: \frac{2 \Dconmax \numobs \lambda_\numobs
 }{\numobs-\spindex-1}, \qquad \mbox{and} \qquad \frac{\sigw}{\Cmax \,
 (\numobs - \spindex-1)} \; \leq \; \Exs[\myuvar_i] \: \leq \:
 \frac{\sigw \Dconmax }{\numobs - \spindex -1}.
\end{equation}
\noindent (b) Moreover, each pair $(\myumean_i, \myuvar_i)$ is sharply
concentrated, in that we have
\begin{eqnarray}
\Prob \Biggr [ |\myumean_i| \geq \frac{6 \Dconmax \numobs
\lambda_n}{\numobs-\spindex-1}, \quad \mbox{or} \quad |\myuvar_i |
\geq 2 \Exs[\myuvar_i] \Biggr ] & \leq & \frac{K}{\numobs - \spindex},
\end{eqnarray}
where $K$ is a fixed constant independent of $\numobs$ and
$\spindex$. \\
\elems
We exploit this lemma as follows.  First define the event
\begin{eqnarray*}
\Tail(\delta) & \defn & \bigcup_{i=1}^\spindex \Biggr \{ |\myumean_i|
\geq \frac{6 \Dconmax \numobs \lambda_n}{\numobs-\spindex-1}, \quad
\mbox{or} \quad |\myuvar_i | \geq 2 \Exs[\myuvar_i] \Biggr \}.
\end{eqnarray*}
By the union bound and Lemma~\ref{LemUbehave}(b), we have
\begin{equation}
\Prob[\Tail(\delta)] \; \leq \; \spindex \frac{K}{\numobs -\spindex}
\nonumber \; = \; \frac{K}{\frac{\numobs}{\spindex} - 1} \;
\rightarrow \; 0,
\end{equation}
since $\frac{\numobs}{\spindex} \rightarrow +\infty$ as $\numobs
\rightarrow +\infty$.  For convenience in notation, for any $a \in
\real$ and $b \in \real_+$, we use $U_i(a, b)$ to denote a Gaussian
random variable with mean $a$ and variance $b$.  Conditioning on the
event $\Tail(\delta)$ and its complement, we have
\begin{eqnarray}
\Prob[\max_{i \in \Sset} \Uvar_i > \rho_\numobs] & \leq &
\Prob[\max_{i \in \Sset} \Uvar_i > \rho_\numobs \; \mid \;
\Tail(\delta)^c ] + \Prob[\Tail(\delta)] \nonumber \\
\label{EqnDoubleTail}
& \leq & \Prob[\max_{i \in \Sset} \Uvar_i(\mustar_i, \vstar_i) >
\rho_\numobs] + \frac{K}{\frac{\numobs}{\spindex} - 1},
\end{eqnarray}
where each $\Uvar_i(\mustar_i, \vstar_i)$ is Gaussian with mean
$\mustar_i \defn 6 \Dconmax \lambda_\numobs
\frac{\numobs}{\numobs-\spindex-1}$ and variance $\vstar_i \defn 2
\Exs[\myuvar_i]$ respectively.  In asserting the
inequality~\eqref{EqnDoubleTail}, we have used the fact that the
probability of the event $\{\max_{i \in \Sset} \myumean_i >
\rho_\numobs\}$ increases as the mean and variance of $\myumean_i$
increase. Continuing the argument, we have
\begin{eqnarray*}
\Prob[\max_{i \in \Sset} \Uvar_i(\mustar_i, \vstar_i) > \rho_\numobs]
& \leq & \Prob[\max_{i \in \Sset} |\Uvar_i(\mustar_i, \vstar_i)| >
\rho_\numobs] \\
& \leq & \frac{1}{\rho_\numobs} \Exs \left[ \max_{i \in \Sset}
|\Uvar_i(\mustar_i, \vstar_i)| \right],
\end{eqnarray*}
where the last step uses Markov's inequality.  We now decompose
$\Uvar_i(\mustar_i, \vstar_i) \stackrel{d}{=} 2 \Dconmax \,
\lambda_\numobs \frac{\numobs}{\numobs-\spindex-1} + \Uvartil_i(0,
\vstar_i)$, and write
\begin{eqnarray*}
\Exs \left[ \max_{i \in \Sset} | \Uvar_i(\mustar_i, \vstar_i)| \right]
& \leq & 2 \Dconmax \, \lambda_\numobs
\frac{\numobs}{\numobs-\spindex-1} + \Exs \left[ \max_{i \in \Sset}
|\Uvartil(0, \vstar_i)| \right].
\end{eqnarray*}
With this decomposition, we use the bound~\eqref{EqnUbehaveBounds} on
$\vstar_i \defn 2 \Exs[\myuvar_i]$ and Lemma~\ref{LemGeneric} on
Gaussian maxima (see Appendix~\ref{AppGaussComp}) to conclude that
\begin{eqnarray*}
\frac{1}{\rho_\numobs} \Exs \left[ \max_{i \in \Sset} |
\Uvar_i(\mustar_i, \vstar_i) | \right] & \leq &
\frac{1}{\rho_\numobs} \left[ 2 \Dconmax \, \lambda_\numobs
\frac{\numobs}{\numobs-\spindex-1} + 3 \sqrt{ \frac{2 \sigw \,
\Dconmax \log \spindex}{\numobs - \spindex -1} } \right],
\end{eqnarray*}
which converges to zero by the second
condition~\eqref{EqnLamGaussCondition} in the theorem statement.

\subsection{Proof of Theorem~\ref{ThmGenNoiseGauss}(a)}

We establish the claim by proving that under the stated conditions,
$\max_{j \in \Sbar} |\Vvar_j |> \lambda_{\numobs}$ with probability
one, for any positive sequence $\lambda_{\numobs} > 0$.  We begin by
writing $\Vvar_j = \Exs[\Vvar_j] + \Vvarcond_j$, where $\Vvarcond_j$
is zero-mean.  Now
\begin{eqnarray*}
\max_{j \in \Sbar} |\Vvar_j| & \geq & \max_{j \in \Sbar} |\Vvarcond_j|
- \max_{j \in \Sbar} |\Exs[\Vvar_j]| \\
& \geq & \max_{j \in \Sbar} |\Vvar_j| - (1 - \matbound)
\lambda_{\numobs}
\end{eqnarray*}
where have used Lemma~\ref{LemGaussCond}.  Consequently, the event $\{
\max_{j \in \Sbar} |\Vvarcond_j| > (2-\matbound) \lambda_\numobs \}$
implies the event $\{\max_{j \in \Sbar} |\Vvar_j| > \lambda_\numobs
\}$, so that
\begin{eqnarray*}
\Prob[\max_{j \in \Sbar} |\Vvar_j| > \lambda_\numobs] & \geq &
\Prob[\max_{j \in \Sbar} |\Vvarcond_j| > (2-\matbound) \,
\lambda_\numobs].
\end{eqnarray*}

From the preceding proof of Theorem~\ref{ThmGenNoiseGauss}(b), we know
that conditioned on $\Amatt{\Sset}$ and $\Esca$, the random vector
$(\Vvar_1, \ldots, \Vvar_\bigN)$ is Gaussian with covariance of the
form $\myVar \, [\CovMat_{\Sbar \Sbar} - \CovMat_{\Sbar
S}(\CovMat_{\Sset \Sset})^{-1} \CovMat_{\Sset \Sbar} ]$; thus, the
zero-mean version $(\Vvarcond_1, \ldots, \Vvarcond_\bigN)$ has the
same covariance.  Moreover, Lemma~\ref{LemExtStats} guarantees that
the random scaling term $\myVar$ is sharply concentrated.  In
particular, defining for any $\delta > 0$ the event $\Tail(\delta)
\defn \{ \, | \myVar - \Exs[\myVar]| \geq \delta \Exs[\myVar] \}$, we
have $\Prob[\Tail(\delta)] \rightarrow 0$, and the bound
\begin{eqnarray*}
\Prob[\max_{j \in \Sbar} |\Vvarcond_j| > (2-\matbound) \,
\lambda_{\numobs}] & \geq & (1-\Prob[\Tail(\delta)]) \; \Prob \left[
\max_{j \in \Sbar} |\Vvarcond_j| > (2-\matbound) \, \lambda_{\numobs}
\; \mid \; \Tail(\delta)^c \right] \\
& \geq & (1 - \Prob[\Tail(\delta)]) \; \Prob \left[ \max_{j \in \Sbar}
  |\Zvar_j(\vstar)| > (2-\matbound) \, \lambda_{\numobs} \right],
\end{eqnarray*}
where each $\Zvar_j \equiv \Zvar_j(\vstar)$ is the conditioned version
of $\Vvarcond_j$ with the scaling factor $\myVar$ fixed to $\vstar
\defn (1-\delta) \Exs[\myVar]$.  (Here we have used the fact that the
probability of Gaussian maxima decreases as the variance decreases,
and that $\var(\Vvarcond_j) \geq \vstar$ when conditioned on
$\Tail(\delta)^c$.)

Our proof proceeds by first analyzing the expected value, and then
exploiting Gaussian concentration of measure.  We summarize the key
results in the following:
\blems
\label{LemExpInfinity}
Under the stated conditions, one of the following two conditions must
hold:
\begin{enumerate}
\item[(a)] either $\frac{\lambda_{\numobs}^2}{\vstar} \rightarrow
+\infty$, and there exists some $\gamma > 0$ such that
$\frac{1}{\lambda_{\numobs}}\Exs[\max_{j \in \Sbar} \Zvar_j] \geq
(2-\matbound) \left[1 + \gamma\right]$ for all sufficiently large
$\numobs$, or
\item[(b)] there exist constants $\alpha, \gamma > 0$ such that
$\frac{\vstar}{\lambda_{\numobs}^2} \leq \alpha$ and
$\frac{1}{\lambda_{\numobs}}\Exs[\max_{j \in \Sbar} \Zvar_j] \geq
\gamma \sqrt{\log \bigN}$ for all sufficiently large $\numobs$.
\end{enumerate}
\elems

\blems
\label{LemGaussConcenLower}
For any $\eta > 0$, we have
\begin{eqnarray}
\label{EqnLipConLower}
\Prob[ \max_{j \in \Sbar} \Zvar_j(\vstar) < \Exs[\max_{j \in \Sbar}
\Zvar_j(\vstar)] - \eta ] & \leq & \exp \left (-\frac{\eta^2}{2
\vstar} \right).
\end{eqnarray}
\elems

Using these two lemmas, we complete the proof as follows.  First, if
condition (a) of Lemma~\ref{LemExpInfinity} holds, then we set $\eta =
\frac{(2-\matbound) \, \gamma \lambda_{\numobs}}{2}$ in
equation~\eqref{EqnLipConLower} to obtain that
\begin{eqnarray*}
\Prob[ \frac{1}{\lambda_{\numobs}} \max_{j \in \Sbar} \Zvar_j(\vstar)
  \geq (2-\matbound) \, (1 + \frac{\gamma}{2}) \, ] & \geq & 1- \exp
  \left (-\frac{(2-\matbound)^2 \, \gamma^2 \lambda_{\numobs}^2 }{8
  \vstar} \right).
\end{eqnarray*}
This probability converges to $1$ since $\frac{\lambda_{\numobs}^2}{\vstar}
\rightarrow +\infty$ from Lemma~\ref{LemExpInfinity}(a).

On the other hand, if condition (b) holds, then we use the bound
$\frac{1}{\lambda_{\numobs}}\Exs[\max_{j \in \Sbar} \Zvar_j] \geq
\gamma \sqrt{\log \bigN}$ and set $\eta = \frac{\gamma\lambda_{\numobs}
\sqrt{\log \bigN}}{2}$ in equation~\eqref{EqnLipConLower} to obtain
\begin{eqnarray*}
\Prob[ \frac{1}{\lambda_{\numobs}} \max_{j \in \Sbar} \Zvar_j(\vstar)
> 2 \, (2-\matbound) \, ] & \geq & \Prob[ \frac{1}{\lambda_{\numobs}}
\max_{j \in \Sbar} \Zvar_j(\vstar) \geq \frac{\gamma \sqrt{\log
\bigN}}{2} \, ]\\
& \geq & 1- \exp \left (-\frac{\gamma^2 \lambda_{\numobs}^2 \log
\bigN}{8 \vstar} \right).
\end{eqnarray*}
This probability also converges to $1$ since
$\frac{\lambda_{\numobs}^2}{\vstar} \geq 1/\alpha$ and $\log \bigN
\rightarrow +\infty$.  Thus, in either case, we have shown that
$\lim_{\numobs \rightarrow +\infty} \Prob[\frac{1}{\lambda_\numobs}
\max_{j \in \Sbar} \Zvar_j(\vstar) > (2-\matbound)] = 1$, thereby
completing the proof of Theorem~\ref{ThmGenNoiseGauss}(a).

\section{Illustrative simulations}
\label{SecExperimental}

In this section, we provide some simulations to confirm the threshold
behavior predicted by Theorem~\ref{ThmGenNoiseGauss}.  We consider the
following three types of sparsity indices:
\begin{enumerate}
\item[(a)] \emph{linear sparsity}, meaning that $\spindex(\mdim) =
\alpha \mdim$ for some $\alpha \in (0,1)$;
\item[(b)] \emph{sublinear sparsity}, meaning that $\spindex(\mdim) =
\alpha \mdim/(\log (\alpha \mdim))$ for some $\alpha \in (0,1)$, and
\item[(c)] \emph{fractional power} sparsity, meaning that
$\spindex(\mdim) =  \alpha \mdim^{\gamma}$ for some $\alpha, \gamma \in (0,1)$.
\end{enumerate}
For all three types of sparsity indices, we investigate the
success/failure of the Lasso in recovering the sparsity pattern, where
the number of observations scales as $\numobs = 2 \, \theta \,
\spindex \log (\mdim - \spindex) + \spindex + 1$.  The \emph{control
parameter} $\theta$ is varied in the interval $(0,2.4)$.  For all
results shown here, we fixed $\alpha = 0.40$ for all three ensembles,
and set $\gamma = 0.75$ for the fractional power ensemble.  In
addition, we set $\lambda_\numobs = \sqrt{ \frac{ \log(\mdim-\spindex)
\log (\spindex)}{\numobs}}$ in all cases.

We begin by considering the uniform Gaussian ensemble, in which each
row $\arow_k$ is chosen in an i.i.d. manner from the multivariate
$N(0, I_{\mdim \times \mdim})$ distribution.  Recall that for the
uniform Gaussian ensemble, the critical value is $\ThreshUp =
\ThreshLow = 1$.  Figure~\ref{FigResultsId} plots the control
parameter $\theta$ versus the probability of success, for linear
sparsity (a), sublinear sparsity pattern (b), and fractional power
sparsity (c), for three different problem sizes ($\mdim \in \{128,
256, 512 \}$). Each point represents the average of $200$ trials.
\newcommand{\mytw}{0.33\textwidth}
\begin{figure}[h]
\begin{center}
\begin{tabular}{ccc}
\widgraph{\mytw}{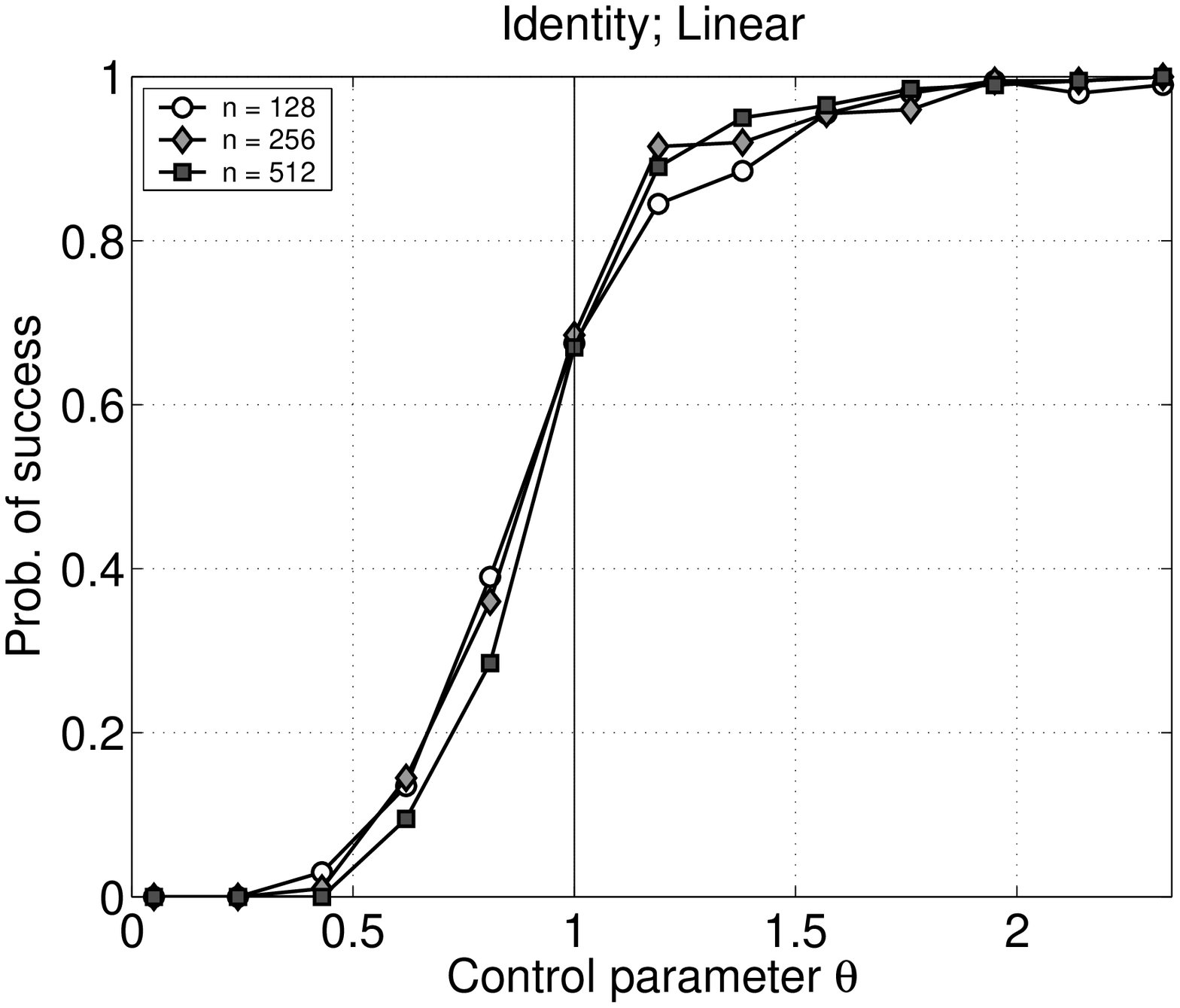} &
\widgraph{\mytw}{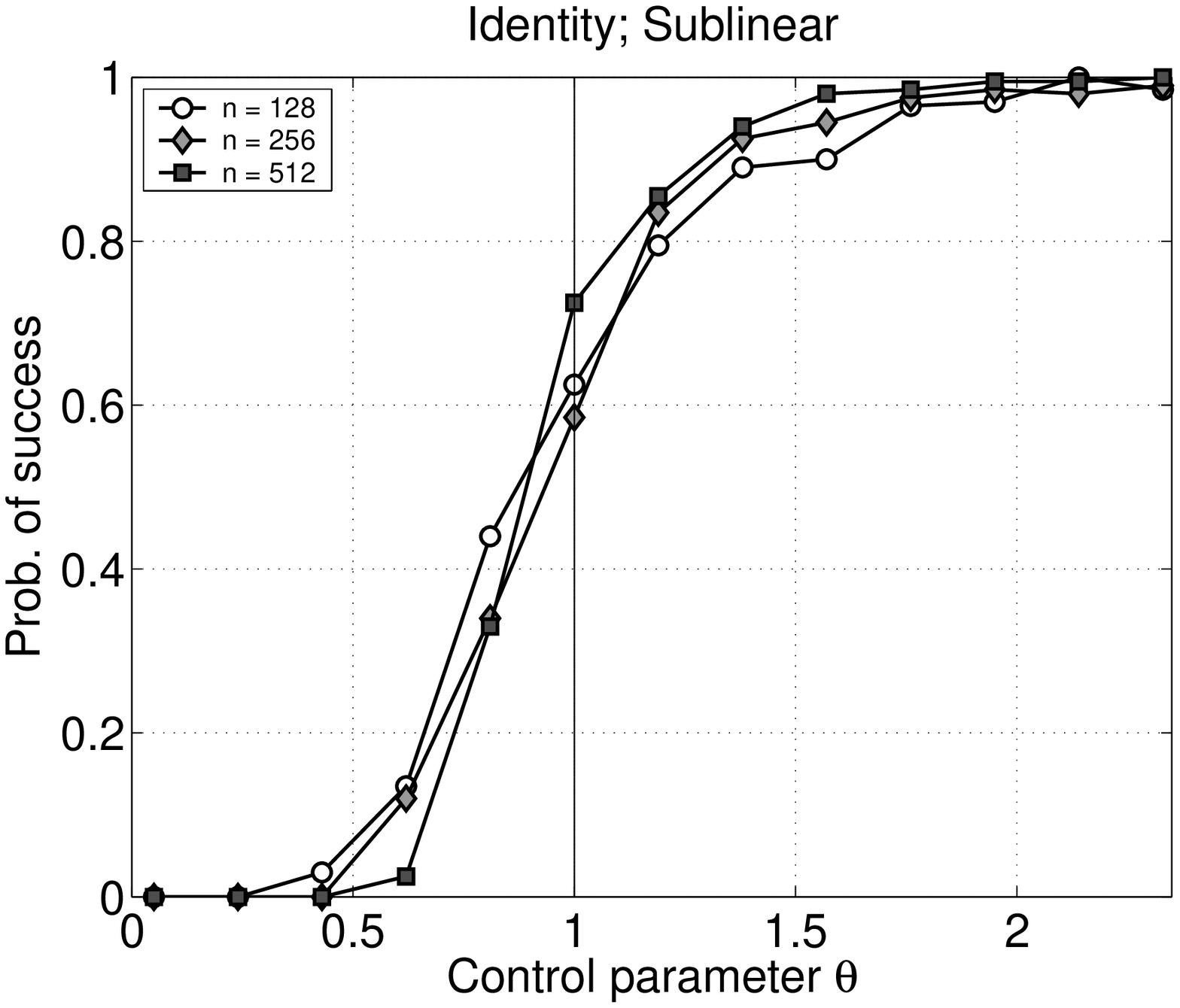} &
\widgraph{\mytw}{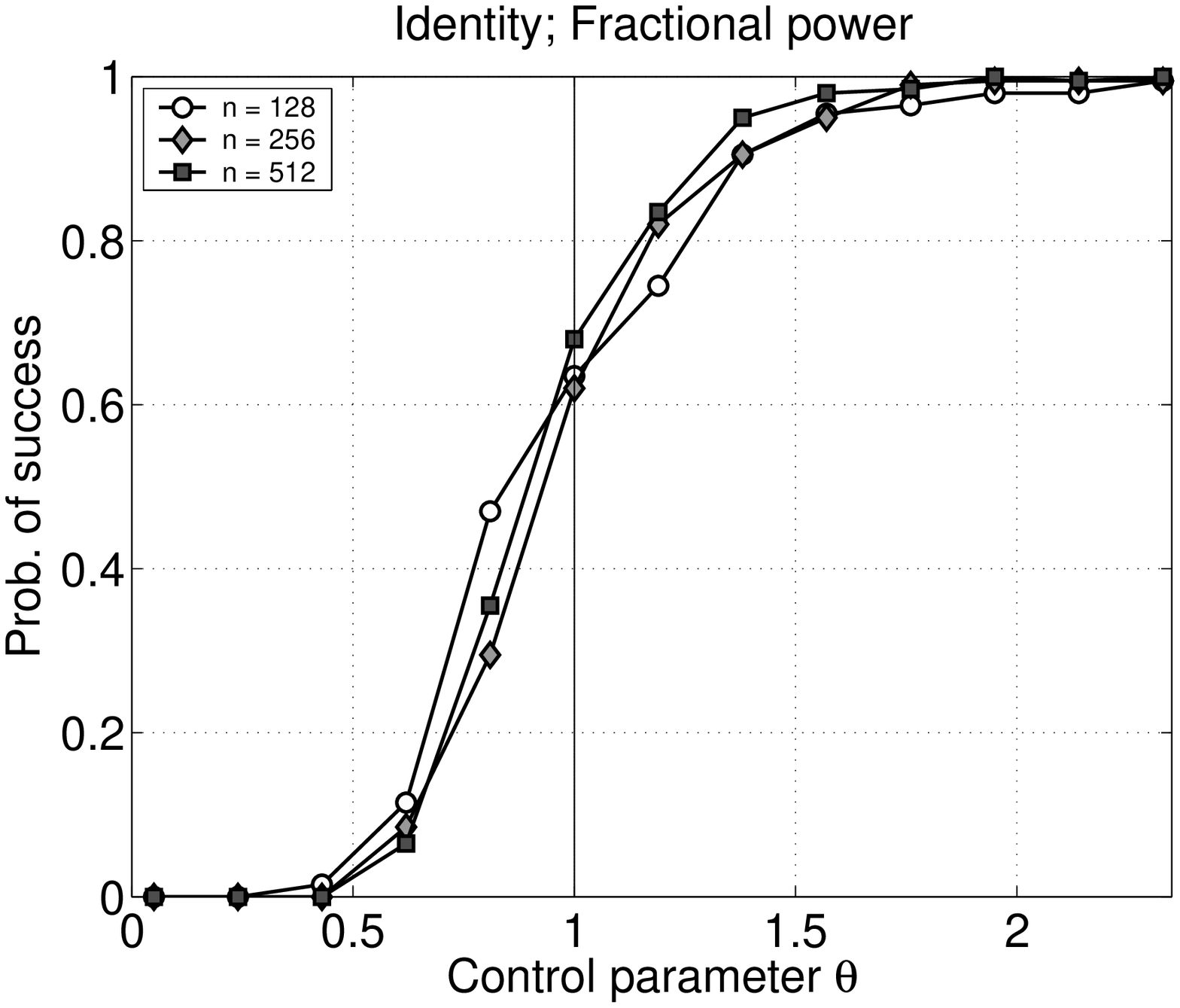} \\
(a) & (b) & (c)
\end{tabular}
\caption{Plots of the number of data samples (indexed by the control
parameter $\theta$ versus the probability of success in the Lasso for
the uniform Gaussian ensemble.  Each panel shows three curves,
corresponding to the problem sizes $\mdim \in \{128, 256, 512 \}$, and
each point on each curve represents the average of $200$ trials.  (a)
Linear sparsity index: $\spindex(\mdim) = \alpha \mdim$.  (b)
Sublinear sparsity index $\spindex(\mdim) = \alpha \mdim/\log(\alpha
\mdim)$.  (c) Fractional power sparsity index $\spindex(\mdim) =
\alpha \mdim^\gamma$ with $\gamma = 0.75$.}
\label{FigResultsId}
\end{center}
\end{figure}
Note how the probability of success rises rapidly from $0$ around the
predicted threshold point $\theta = 1$, with the sharpness of the
threshold increasing for larger problem sizes.

We now consider a non-uniform Gaussian ensemble---in particular, one
in which the covariance matrices $\CovMat$ are Toeplitz with the
structure
\begin{eqnarray}
\label{EqnToeplitzFamily}
\CovMat & = & \begin{bmatrix} 1 & \toep & \toep^2 & \cdots &
\toep^{\mdim-1} & \toep^\mdim \\
\toep & 1 & \toep & \toep^2 & \cdots &  \toep^{\mdim-1} \\
\toep^2 & \toep & 1 & \toep & \cdots & \toep^{\mdim-2} \\
\vdots & \vdots & \vdots & \vdots & \vdots &  \vdots \\
\toep^{\mdim} & \cdots & \toep^3 & \toep^2 & \toep & 1
	      \end{bmatrix},
\end{eqnarray}
for some $\toep \in (-1, +1)$.  As shown by Zhao and Yu~\cite{Zhao06},
this family of Toeplitz matrices satisfy condition~\eqref{EqnCovCond}.
Moreover, the maximum and minimum eigenvalues ($\Cmin$ and $\Cmax$)
can be computed using standard asymptotic results on Toeplitz matrix
families~\cite{GrayToeplitz}.
\begin{figure}[h]
\begin{center}
\begin{tabular}{ccc}
\widgraph{\mytw}{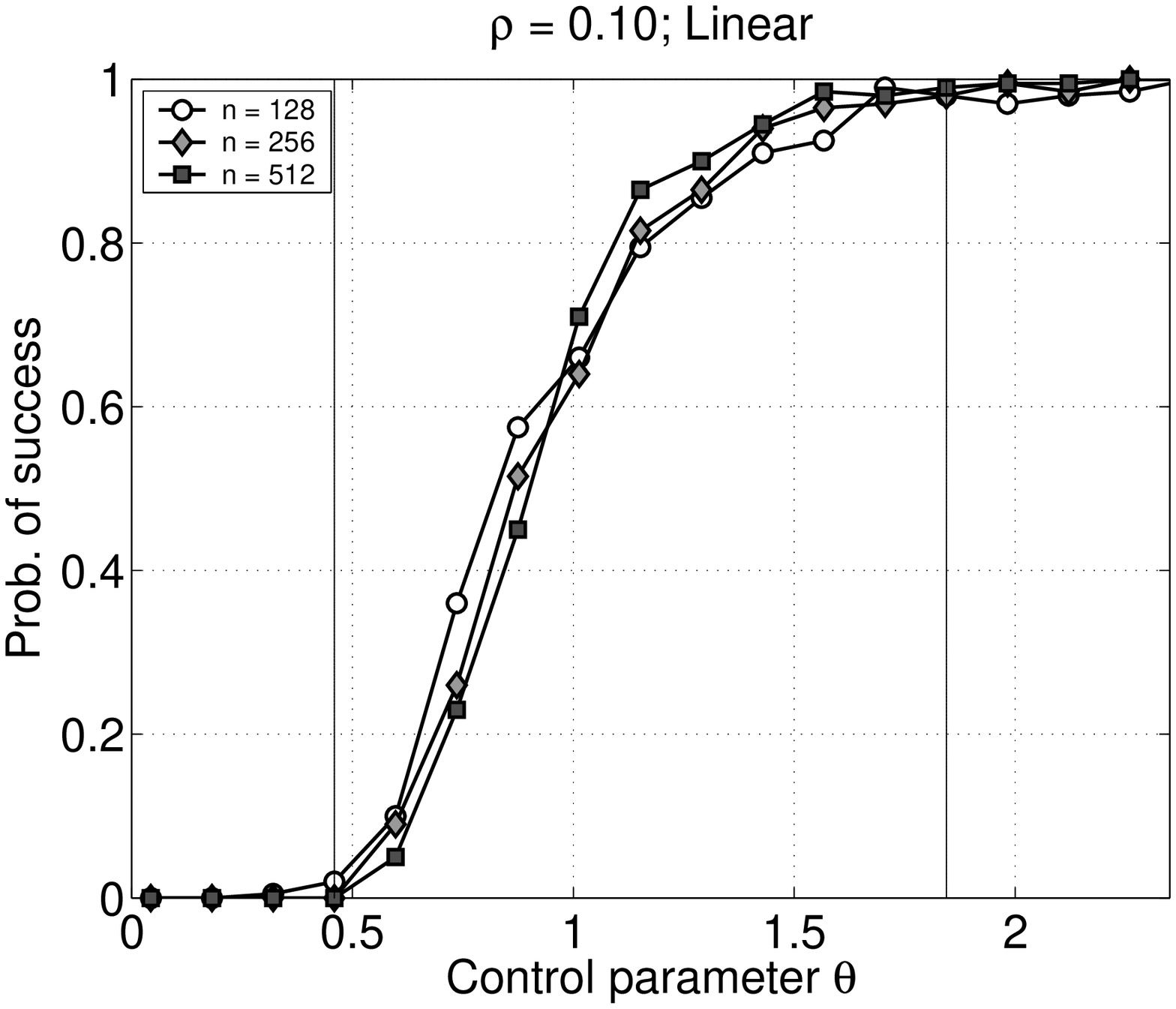} &
\widgraph{\mytw}{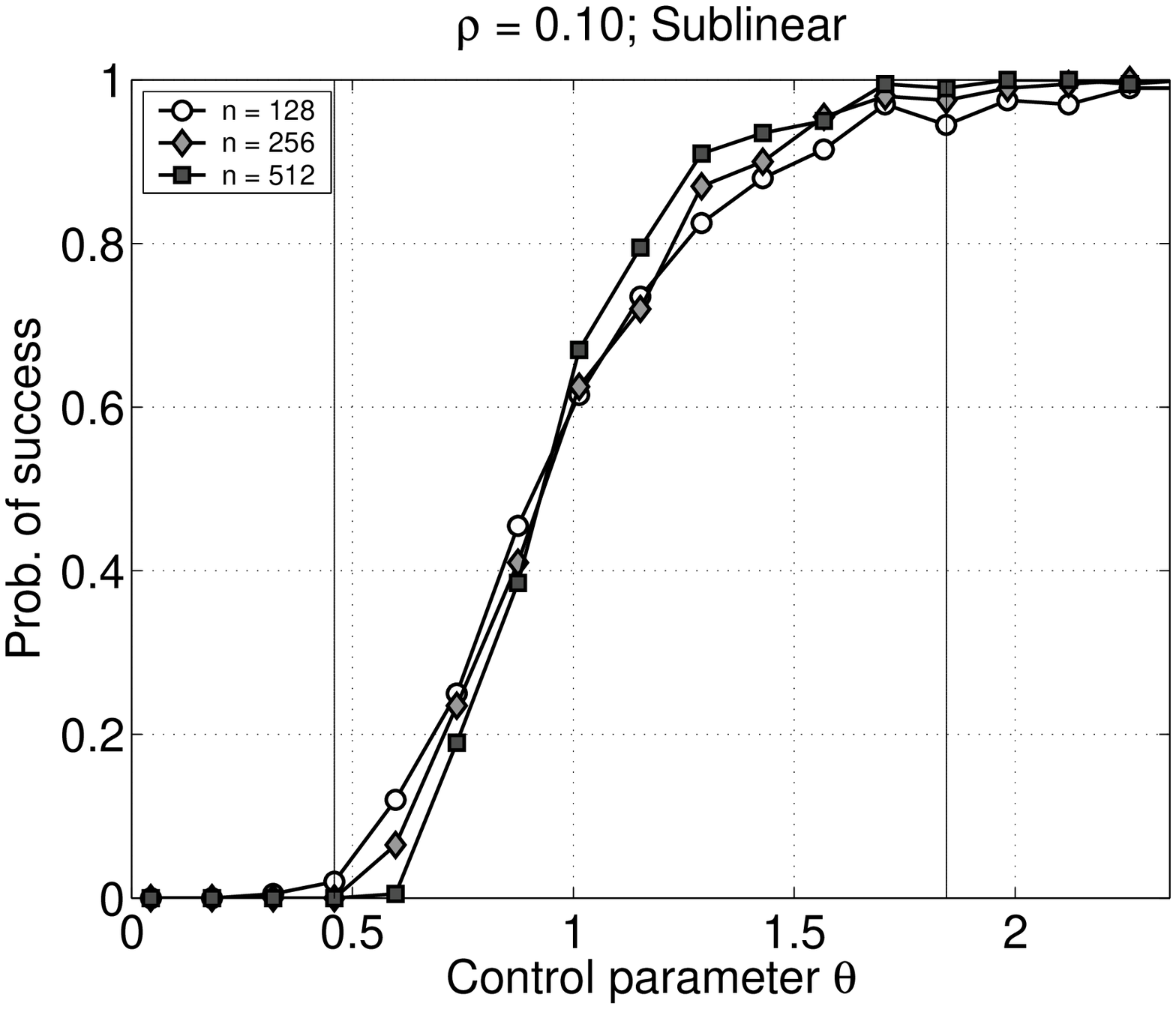} &
\widgraph{\mytw}{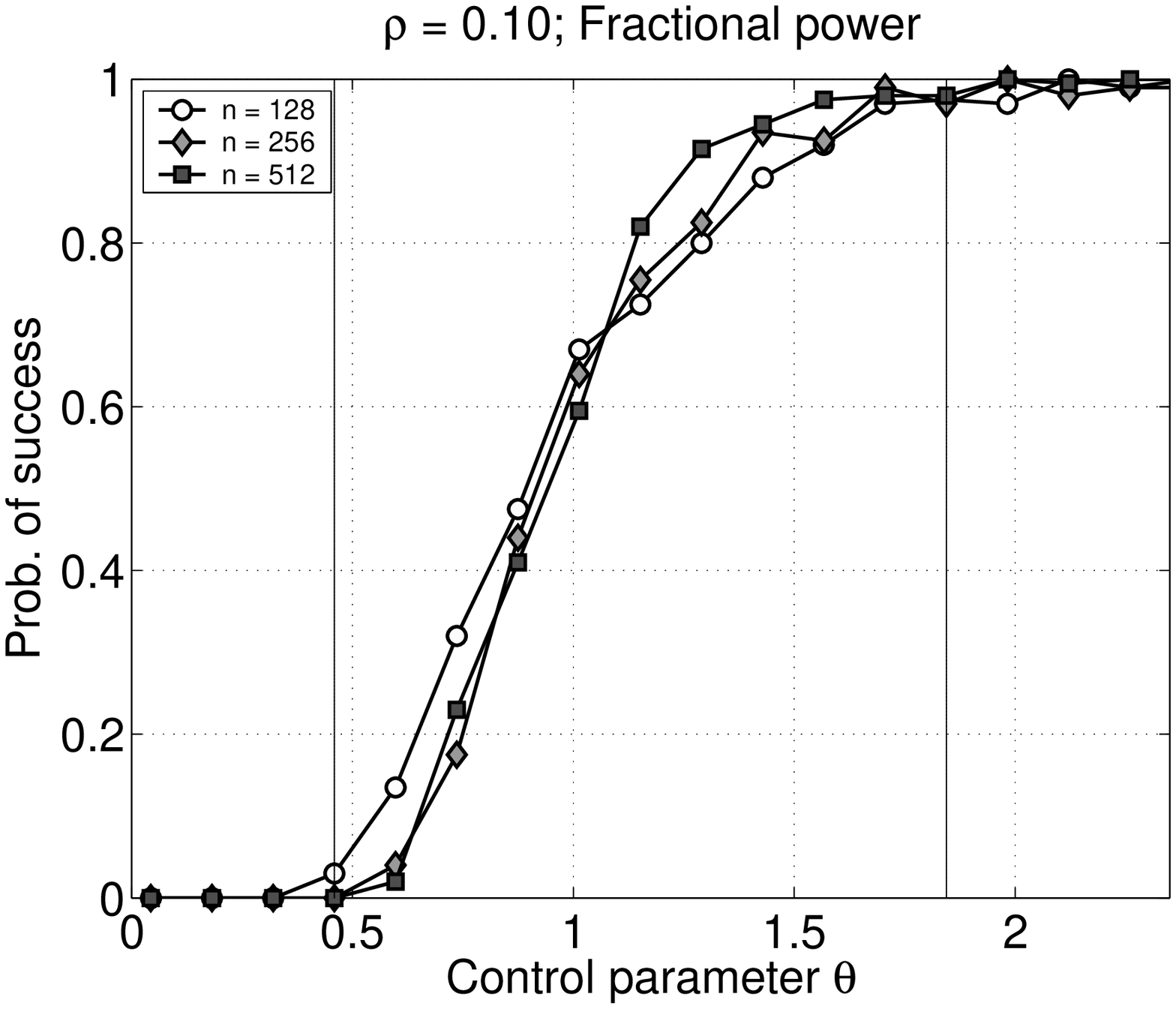} \\
(a) & (b) & (c)
\end{tabular}
\caption{Plots of the number of data samples (indexed by the control
parameter $\theta$ versus the probability of success in the Lasso for
the Toeplitz family~\eqref{EqnToeplitzFamily} with $\rho = 0.10$.
Each panel shows three curves, corresponding to the problem sizes
$\mdim \in \{128, 256, 512 \}$, and each point on each curve
represents the average of $200$ trials.  (a) Linear sparsity index:
$\spindex(\mdim) = \alpha \mdim$.  (b) Sublinear sparsity index
$\spindex(\mdim) = \alpha \mdim/\log(\alpha \mdim)$. (c) Fractional
power sparsity index $\spindex(\mdim) = \alpha \mdim^\gamma$ with
$\gamma = 0.75$.}
\label{FigResultsRho}
\end{center}
\end{figure}
Figure~\ref{FigResultsRho} shows representative results for this
Toeplitz family with $\rho = 0.10$.  Panel (a) corresponds to linear
sparsity $\spindex = \alpha \mdim$ with $\alpha = 0.40$), and panel
(b) corresponds to sublinear sparsity ($\spindex = \alpha
\mdim/\log(\alpha \mdim)$ with $\alpha = 0.40$).  Each panel shows
three curves, corresponding to the problem sizes $\mdim \in \{128,
256, 512 \}$, and each point on each curve represents the average of
$200$ trials.  The vertical lines to the left and right of $\theta =
1$ represent the theoretical upper and lower bounds on the threshold
($\ThreshUp \approx 1.84$ and $\ThreshLow \approx 0.46$ respectively
in this case).  Once again, these simulations show good agreement with
the theoretical predictions.

\section{Discussion}
\label{SecDiscussion}

The problem of recovering the sparsity pattern of a high-dimensional
vector $\betastar$ from noisy observations has important applications
in signal denoising, graphical model selection, sparse approximation,
and subset selection.  This paper focuses on the behavior of
$\ell_1$-regularized quadratic programming, also known as the Lasso,
for estimating such sparsity patterns in the noisy and
high-dimensional setting. The main contribution of this paper is to
establish a set of general and sharp conditions on the observations
$\numobs$, the sparsity index $\spindex$ (i.e., number of non-zero
entries in $\betastar$), and the ambient dimension $\mdim$ that
characterize the success/failure behavior of the Lasso in the
high-dimensional setting, in which $\numobs$, $\mdim$ and $\spindex$
all tend to infinity.  For the uniform Gaussian ensemble, our
threshold result is sharp, whereas for more general Gaussian
ensembles, it should be possible to tighten the analysis given here.

\subsection*{Acknowledgements}
We would like to thank Noureddine El Karoui and Bin Yu for helpful comments and
pointers.  This work was partially supported by an Alfred P. Sloan
Foundation Fellowship, and an Intel Corporation Equipment Grant.


\appendix

\section{Proof of Lemma~\ref{LemNecSuffCond}}
\label{AppNecSuff}
By standard conditions for optimality in a convex
program~\cite{Hiriart1}, the point $\betahat \in \real^\mdim$ is
optimal if and only if there exists a subgradient $\subgradopt \in
\partial \ell_1(\betahat)$ such that
\begin{eqnarray}
\label{EqnOpt1}
\frac{1}{\numobs} \Amat^T \Amat \betahat - \frac{1}{\numobs} \Amat^T
\ysca + \lambda \subgradopt & = & 0.
\end{eqnarray}
Here the subdifferential of the $\ell_1$ norm takes the form
\begin{equation*}
\partial \ell_1(\betahat) = \left \{ \subgradopt \in \real^\mdim \;
\mid \; \subgradopt_i = \sgn(\betahat_i) \quad \mbox{for $\betahat_i
\neq 0$}, \qquad |\subgradopt_j| \leq 1 \quad \mbox{otherwise} \right
\}.
\end{equation*}
Substituting our observation model $\ysca = \Amat \betastar + \esca$
and re-arranging yields
\begin{eqnarray}
\label{EqnKeyInter}
\frac{1}{\numobs} \Amat^T \Amat (\betahat - \betastar) -
\frac{1}{\numobs} \Amat^T \esca + \lambda \subgradopt & = & 0.
\end{eqnarray}
Now condition $\Prop(\Amat, \betastar, \wsca, \lambda)$ holds if and
only we have
\begin{equation*}
\betahat_{\Sbar} = 0, \quad \betahat_{\Sset} \neq 0, \qquad \mbox{and}
\qquad \subgradopt_\Sset = \sgn(\betastar_\Sset), \quad \quad
|\subgradopt_{\Sbar}| \leq 1.
\end{equation*}
From these conditions and using equation~\eqref{EqnKeyInter}, we
conclude that the condition $\Prop(\Amat, \betastar,\wsca, \lambda)$
holds if and only if
\begin{subequations}
\begin{eqnarray*}
\frac{1}{\numobs} \Amatt{\Sbar}^T \Amatt{\Sset} \left(\betahat_\Sset -
\betastar_\Sset \right) - \frac{1}{n} \Amatt{\Sbar}^T \esca & = & -
\lambda \subgradopt_{\Sbar}. \\
\frac{1}{\numobs} \Amatt{\Sset}^T \Amatt{\Sset} \left(\betahat_\Sset - \betastar_\Sset
\right) - \frac{1}{\numobs}\Amatt{\Sset}^T \esca & = & - \lambda
\sgn(\betastar_\Sset).
\end{eqnarray*}
\end{subequations}
Using the invertibility of $\Amatt{\Sset}^T \Amatt{\Sset}$, we may solve for
$\betahat_\Sset$ and $\subgradopt_{\Sbar}$ to conclude that
\begin{subequations}
\begin{eqnarray*}
\lambda \; \subgradopt_{\Sbar} & = & \Amatt{\Sbar}^T \Amatt{\Sset}
\left(\Amatt{\Sset}^T \Amatt{\Sset} \right)^{-1}
\left[\frac{1}{\numobs} \Amatt{\Sset}^T \esca- \lambda
\sgn(\betastar_\Sset) \right] - \frac{1}{\numobs} \Amatt{\Sbar}^T
\esca \\
\betahat_\Sset & = & \betastar_\Sset + \left(\frac{1}{\numobs}
\Amatt{\Sset}^T \Amatt{\Sset} \right)^{-1} \left[\frac{1}{\numobs}
\Amatt{\Sset}^T \esca - \lambda \sgn(\betastar_\Sset) \right].
\end{eqnarray*}
\end{subequations}
From these relations, the conditions $| \subgradopt_{\Sbar} | \leq 1$
and $\betahat_\Sset \neq 0$ yield conditions~\eqref{EqnPropA}
and~\eqref{EqnPropB} respectively.

\section{Some Gaussian comparison results}
\label{AppGaussComp}
We state here (without proof) some well-known comparison results on
Gaussian maxima~\cite{LedTal91}.  We begin with a crude but useful
bound:
\blems
\label{LemGeneric}
For any Gaussian random vector $(X_1, \ldots, X_n)$, we have
\begin{eqnarray*}
\Exs \max_{1 \leq i \leq n} |X_i| & \leq & 3 \sqrt{\log n} \; \max_{1
\leq i \leq n} \sqrt{\Exs X_i^2}.
\end{eqnarray*}
\elems
\noindent Next we state (a version of) the Sudakov-Fernique
inequality~\cite{LedTal91,Chatterjee05}:
\blems
\label{LemSudFer}
Let $X = (X_1, \ldots, X_n)$ and $Y = (Y_1, \ldots, Y_n)$ be Gaussian
random vectors such that for all $i, j$
\begin{eqnarray*}
\Exs[(Y_i - Y_j)^2] & \leq & \Exs[(X_i - X_j)^2].
\end{eqnarray*}
Then $\Exs[\max \limits_{1 \leq i \leq n} Y_i] \leq \Exs[\max
\limits_{1 \leq i \leq n} X_i]$.
\elems

\section{Auxiliary lemma}
\label{AppAux}

For future use, we state formally the following elementary
\blems
\label{LemSimpleInequal}
Given a collection $\{\Zvar_1, \Zvar_2, \ldots, \Zvar_\bigN\}$ of
zero-mean random variables, for any constant $a > 0$ we have
\begin{subequations}
\begin{eqnarray}
\label{EqnBoundA}
\Prob [\max_{1 \leq j \leq \bigN} |\Zvar_j| \leq a] & \leq & \Prob[\max_{1
\leq j \leq \bigN} \Zvar_j \leq a], \qquad \mbox{and} \\
\label{EqnBoundB}
\Prob [\max_{1 \leq j \leq \bigN} |\Zvar_j| > a] & \leq &
2 \Prob[\max_{1 \leq j \leq \bigN} \Zvar_j > a].
\end{eqnarray}
\end{subequations}
\spro
The first inequality is trivial.  To establish the
inequality~\eqref{EqnBoundB}, we write
\begin{eqnarray*}
\Prob [\max_{1 \leq j \leq \bigN} |\Zvar_j| > a] & = & \Prob [(\max_{1 \leq j
\leq \bigN} \Zvar_j > a) \; \mbox{or} \; (\min_{1 \leq j \leq \bigN} \Zvar_j < -a)] \\
& \leq & \Prob [\max_{1 \leq j \leq \bigN} \Zvar_j > a] + \Prob[\min_{1 \leq j
\leq \bigN} \Zvar_j < -a] \\
& = & 2 \Prob [\max_{1 \leq j \leq \bigN} \Zvar_j > a],
\end{eqnarray*}
where we have used the union bound, and the symmetry of the events
$\{\max_{1 \leq j \leq \bigN} \Zvar_j > a\}$ and $\{\min_{1 \leq j \leq \bigN}
\Zvar_j < -a \}$.
\fpro
\elems

\section{Lemma for Theorem~\ref{ThmGenNoiseGauss}}

\label{AppOneTwo}

\subsection{Proof of Lemma~\ref{LemGaussCond}}
\label{AppGaussCond}

Conditioned on both $\Amatt{\Sset}$ and $\Esca$, the only random
component in $\Vvar_j$ is the column vector $\acol_j$.  Using standard
LLSE formula~\cite[e.g.,]{Bickel01} (i.e., for estimating
$\Amatt{\Sbar}$ on the basis of $\Amatt{\Sset}$), the random variable
$(\Amatt{\Sbar} \, \mid \, \Amatt{\Sset}, \Esca) \sim (\Amatt{\Sbar}\, \mid \, \Amatt{\Sset})$ is Gaussian with mean and covariance
\begin{subequations}
\begin{eqnarray}
\Exs[\Amatt{\Sbar}^T \; \mid \; \Amatt{\Sset}, \; \Esca] & = &
\CovMat_{\Sbar \Sset} (\CovMat_{\Sset\Sset})^{-1} \Amatt{\Sset}^T, \\
\var(\Amatt{\Sbar} \, \mid \, \Amatt{\Sset}) & = & \CondMat \; = \;
\CovMat_{\Sbar \Sbar} - \CovMat_{\Sbar \Sset} (\CovMat_{\Sset\Sset})^{-1}
\CovMat_{\Sset \Sbar}.
\end{eqnarray}
\end{subequations}
Consequently, we have
\begin{eqnarray*}
\left | \Exs[V_j \, \mid \, \Amatt{\Sset}, W] \right| & = & \left |
\CovMat_{\Sbar \Sset} (\CovMat_{\Sset\Sset})^{-1} \Amatt{\Sset}^T
\Biggr \{ \Amatt{\Sset} \left(\Amatt{\Sset}^T \Amatt{\Sset}
\right)^{-1} \lambda_{\numobs} \spec -
\left[\Amatt{\Sset}\left(\Amatt{\Sset}^T \Amatt{\Sset} \right)^{-1}
\Amatt{\Sset}^T - I_{\numobs \times \numobs} \right]
\frac{\Esca}{\numobs} \Biggr \} \right | \nonumber \\
& = & \left | \CovMat_{\Sbar \Sset} (\CovMat_{\Sset \Sset})^{-1}
\lambda_{\numobs} \spec \right| \nonumber \\
& \leq & \lambda_{\numobs} (1-\matbound) \ones,
\end{eqnarray*}
as claimed.

Similarly, we compute the elements of the conditional covariance
matrix as follows
\begin{multline*}
\cov(V_j, V_k \, \big | \, \Amatt{\Sset}, \Esca) = \\
 \cov(A_{ji}, A_{ki} \, \mid \, \Amatt{\Sset}, \Esca) \; \left \{
\lambda_{\numobs}^2 \spec^T (\Amatt{\Sset}^T \Amatt{\Sset})^{-1} \spec
+ \frac{1}{n^2} \Esca^T \left[I_{\numobs \times \numobs} -
\Amatt{\Sset} \left(\Amatt{\Sset}^T \Amatt{\Sset} \right)^{-1}
\Amatt{\Sset}^T \right] \Esca \right \}.
\end{multline*}

\subsection{Proof of Lemma~\ref{LemExtStats}}
\label{AppExtStats}

We begin by computing the expected value.  Since $\Amatt{\Sset}^T
\Amatt{\Sset}$ is Wishart with matrix $\CovMat_{\Sset\Sset}$, the
random matrix $(\Amatt{\Sset}^T \Amatt{\Sset})^{-1}$ is inverse
Wishart with mean $\Exs[(\Amatt{\Sset}^T \Amatt{\Sset})^{-1}] = \frac
{(\CovMat_{\Sset\Sset})^{-1} }{n-s -1}$ (see Lemma 7.7.1 of
Anderson~\cite{AndersonStat}).  Hence we have
\begin{equation}
\label{EqnInvChi}
\Exs \left[\lambda_{\numobs}^2 \spec^T \left(\Amatt{\Sset}^T
\Amatt{\Sset} \right)^{-1} \spec \right] = \frac{\lambda_{\numobs}^2
}{n-s-1} \; \spec^T (\CovMat_{\Sset\Sset})^{-1} \spec.
\end{equation}
Now define the random matrix $R = I_{\numobs \times \numobs} -
\Amatt{\Sset} (\Amatt{\Sset}^T \Amatt{\Sset})^{-1} \Amatt{\Sset}^T$.
A straightforward calculation yields that $R^2 = R$, so that all the
eigenvalues of $R$ are either $0$ or $1$.  In particular, for any
vector $z = \Amatt{\Sset} u$ in the range of $\Amatt{\Sset}$, we have
\begin{equation}
R z = \left[ I_{\numobs \times \numobs} - \Amatt{\Sset}
(\Amatt{\Sset}^T \Amatt{\Sset})^{-1} \Amatt{\Sset}^T \right ]
\Amatt{\Sset} u \; = \; 0.
\end{equation}
Hence $\dim(\ker R) = \dim(\range \Amatt{\Sset}) = s$.  Since $R$ is
symmetric and positive semidefinite, there exists an orthogonal matrix
$U$ such that $R = U^T D U$, where $D$ is diagonal with $(n-s)$ ones,
and $s$ zeros.  The random matrices $D$ and $U$ are both independent
of $\Esca$, since $\Amatt{\Sset}$ is independent of $\Esca$.  Hence we
have
\begin{eqnarray}
\frac{1}{\numobs^2} \Exs \left[ \Esca^T R \Esca \; \mid \;
\Amatt{\Sset} \right] & = & \frac{1}{\numobs^2} \Exs \left[ \Esca^T
U^T D U \Esca \; \mid \Amatt{\Sset} \right] \nonumber \\
\label{EqnRexp}
& = & \frac{1}{\numobs^2} \trace D U U^T \Exs\left[ \Esca \Esca^T \mid
 \; \Amatt{\Sset} \right ] \nonumber \\ 
\label{EqnCalcWRW}
& = & \sigw \, \frac{\numobs-\spindex}{\numobs^2}
\end{eqnarray}
since $\Exs[\Esca \Esca^T] = \sigw I$.  Consequently, we have
established that $\Exs[\myVar] =
\frac{\lambda_{\numobs}^2}{\numobs-\spindex-1} \; \spec^T
(\CovMat_{\Sset\Sset})^{-1} \spec + \frac{\sigw \,
(\numobs-\spindex)}{\numobs^2}$ as claimed.

We now compute the expected value of the squared variance
\begin{eqnarray*}
\myVar^2 & = & \underbrace{\lambda_{\numobs}^4 \left[ \spec^T
\left(\Amatt{\Sset}^T \Amatt{\Sset} \right)^{-1} \spec \right]^2} +
\underbrace{2 \frac{\lambda_{\numobs}^2}{\numobs^2} \left[\spec^T
\left( \Amatt{\Sset}^T \Amatt{\Sset} \right)^{-1} \spec \right] \left(
\Esca^T R \Esca \right)} + \underbrace{\frac{1}{\numobs^4}
\left(\Esca^T R \Esca \right)^2 } \\
& & \qquad \qquad \myTemp_1 \qquad \qquad \qquad \qquad \qquad
\myTemp_2 \qquad \qquad \qquad \qquad \qquad \qquad \myTemp_3
\end{eqnarray*}
First, conditioning on $\Amatt{\Sset}$ and using the eigenvalue
decomposition $D$ of $R$, we have
\begin{eqnarray}
\Exs[\myTemp_3 | \Amatt{\Sset}] & = & \frac{1}{\numobs^4}
\Exs[(\Esca^T D \Esca)^2] \nonumber \\
& = & \frac{1}{\numobs^4} \Exs\left[(\sum_{i=1}^{\numobs-\spindex}
\Esca_i)^2\right] \nonumber \\
\label{EqnM3}
& = & \frac{2 (\numobs - \spindex) \sigwfour}{\numobs^4} +
\frac{(\numobs - \spindex)^2 \sigwfour}{\numobs^4}.
\end{eqnarray}
whence $\Exs[\myTemp_3] = \frac{2 (\numobs - \spindex)
\sigwfour}{\numobs^4} + \frac{(\numobs - \spindex)^2
\sigwfour}{\numobs^4}$ as well.

Similarly, using conditional expectation and our previous
calculation~\eqref{EqnCalcWRW} of $\Exs[\Esca^T R \Esca \, \mid \,
\Amatt{\Sset}]$, we have
\begin{eqnarray}
\label{EqnM2}
\Exs[\myTemp_2] & = & \frac{2 \lambda_\numobs^2}{\numobs^2} \Exs
\Biggr[ \Exs \Big[\spec^T (\Amatt{\Sset}^T \Amatt{\Sset})^{-1} \spec \;
(\Esca^T R \Esca) \; \mid \; \Amatt{\Sset} \Big] \Biggr] \nonumber \\
& = & \frac{2 \lambda_\numobs^2 \, (\numobs - \spindex)
\sigw}{\numobs^2} \Exs \Big[\spec^T (\Amatt{\Sset}^T \Amatt{\Sset})^{-1} \spec
\Big] \nonumber \\
& = & \frac{2 \lambda_{\numobs}^2 \, (\numobs - \spindex) \,
\sigw}{\numobs^2 \, (\numobs-\spindex-1)} \spec^T
\left(\CovMat_{\Sset\Sset} \right)^{-1} \spec,
\end{eqnarray}
where the final step uses Lemma 7.7.1 of Anderson~\cite{AndersonStat}
on the expectation of inverse Wishart matrices.

Lastly, since $(\Amatt{\Sset}^T \Amatt{\Sset})^{-1}$ is inverse
Wishart with matrix $(\CovMat_{\Sset\Sset})^{-1}$, we can use formula
for second moments of inverse Wishart matrices (see, e.g.,
Siskind~\cite{Siskind72}) to write, for all $\numobs > \spindex +3$,
\begin{eqnarray*}
\Exs[\myTemp_1] & = & \frac{\lambda_{\numobs}^4}{(\numobs-\spindex)\,
 (\numobs-\spindex-3)} \; \left[\spec^T (\CovMat_{\Sset\Sset})^{-1}
 \spec\right]^2 \left \{ 1 + \frac{1}{\numobs-\spindex-1} \right \}.
\end{eqnarray*}
Consequently, combining our results, we have
\begin{eqnarray}
\var(\myVar) & = & \Exs[\myVar^2] - (\Exs[\myVar])^2 \nonumber \\
& = & \sum_{i=1}^3 \Exs[\myTemp_i] - \left\{ \frac{\sigwfour (\numobs
  - \spindex)^2}{\numobs^4} + 2 \frac{\sigw
  (\numobs-\spindex)}{\numobs^2} \frac{
  \lambda_{\numobs}^2}{\numobs-\spindex-1} \; \spec^T
(\CovMat_{\Sset\Sset})^{-1} \spec +
\left(\frac{\lambda_{\numobs}^2}{\numobs-\spindex-1} \; \spec^T
(\CovMat_{\Sset\Sset})^{-1} \spec \right)^2 \right\} \nonumber \\
\label{EqnHvarr}
& = & \underbrace{\frac{2 (\numobs-\spindex) \sigwfour}{\numobs^4}} +
\underbrace{\frac{\lambda_{\numobs}^4 \; [\spec^T
(\CovMat_{\Sset\Sset})^{-1} \spec ]^2}{(\numobs-\spindex-1) \,
(\numobs-\spindex-3) } \left \{ \frac{1}{(\numobs-\spindex)} +
\frac{\numobs-\spindex-1}{(\numobs-\spindex)} - \frac{(\numobs -
\spindex -3)}{(\numobs-\spindex-1)} \right \} }. \\
& & \qquad \myVarr_1 \qquad \qquad \qquad \qquad \qquad \myVarr_2 \nonumber
\end{eqnarray}

Finally, we establish the concentration result.  Using Chebyshev's inequality,
we have
\begin{eqnarray*}
\Prob \left[ |\myVar - \Exs[\myVar]| \geq \delta \Exs[\myVar] \right]
& \leq & \frac{\var(\myVar)}{\delta^2 (\Exs[\myVar])^2}, 
\end{eqnarray*}
so that it suffices to prove that $\var(\myVar)/(\Exs[\myVar])^2
\rightarrow 0$ as $n \rightarrow +\infty$.  We deal with each of the
two variance terms $\myVarr_1$ and $\myVarr_2$ in
equation~\eqref{EqnHvarr} separately.  First, we have
\begin{eqnarray*}
\frac{\myVarr_1}{(\Exs[\myVar])^2} & \leq & \frac{2 (\numobs-\spindex)
\sigwfour}{n^4} \frac{n^4}{(\numobs-\spindex)^2 \sigwfour} \; = \frac{2}{\numobs-\spindex} \;
\rightarrow \; 0.
\end{eqnarray*}
Secondly, denoting $A = (\spec^T (\Amatt{\Sset}^T \Amatt{\Sset})^{-1}
\spec)$ for short-hand, we have
\begin{eqnarray*}
\frac{\myVarr_2}{(\Exs[\myVar])^2} & \leq &
\frac{(\numobs-\spindex-1)^2}{\lambda_{\numobs}^4 A^2}
\frac{\lambda_{\numobs}^4 A^2}{(\numobs-\spindex-1) \,
(\numobs-\spindex-3)} \; \left \{ \frac{1}{(\numobs-\spindex)} +
\frac{\numobs-\spindex-1}{(\numobs-\spindex)} -
\frac{(\numobs-\spindex-3)}{(\numobs-\spindex-1)} \right \} \\
& = & \frac{(\numobs-\spindex-1)}{(\numobs-\spindex-3)} \left \{ \frac{1}{(\numobs-\spindex)} + \frac{\numobs-\spindex-1}{(\numobs-\spindex)} 
 - \frac{(\numobs-\spindex-3)}{(\numobs-\spindex-1)} \right \},
\end{eqnarray*}
which also converges to $0$ as $(\numobs-\spindex) \rightarrow 0$.


\subsection{Proof of Lemma~\ref{LemOne}}

Recall that the Gaussian random vector $(Z_1, \ldots, Z_\bigN)$ is
zero-mean with covariance $\vstar \CondMat$, where $\CondMat \defn
\CovMat_{\Sbar\Sbar} - \CovMat_{\Sbar \Sset}(\CovMat_{\Sset
\Sset})^{-1} \CovMat_{\Sset \Sbar}$. For any index $i$, let $e_i \in
\real^\bigN$ be equal to $1$ in position $i$, and zero otherwise.  For
any two indices $i \neq j$, we have
\begin{eqnarray*}
\Exs[(\Zvar_i - \Zvar_j)^2] & = & \vstar (e_i - e_j)^T \CondMat (e_i -
e_j) \\
& \leq & 2 \vstar \lambda_{max}(\CondMat) \\
& \leq & 2 \Cmax \vstar,
\end{eqnarray*}
since $\CondMat \preceq \CovMat_{\Sbar \Sbar}$ by definition, and
$\myeigmax(\CovMat_{\Sbar \Sbar}) \leq \myeigmax(\CovMat) \leq \Cmax$.

Letting $(X_1, \ldots, X_\bigN) \sim N(0, \Cmax \vstar I_{\bigN \times
\bigN})$, we have $\Exs[(X_i - X_j)^2] = 2 \Cmax \vstar$.  Hence,
applying the Sudakov-Fernique inequality~\cite{LedTal91} yields
$\Exs[\max_j \Zvar_j] \leq \Exs[ \max_j X_j ]$.  By asymptotic
behavior of i.i.d. Gaussians~\cite{GalambosExt, DavNag03}, we have
$\lim_{\bigN \rightarrow \infty} \frac{\Exs[\max_{j} X_j]}{\sqrt {2
\Cmax \vstar \log \bigN}} = 1$.  Consequently, for all $\interdelta >
0$, there exists an $\bigN(\interdelta)$ such that for all $\bigN \geq
\bigN(\interdelta)$, we have
\begin{eqnarray*}
\frac{1}{\lambda_{\numobs}} \Exs[\max_j \Zvar_j(\vstar)] & \leq &
\frac{1}{\lambda_{\numobs}} \Exs[\max_j X_j] \\
& \leq & (1+\interdelta) \; \sqrt{\frac{2 \Cmax \vstar \log
\bigN}{\lambda_{\numobs}^2}} \\
& = & (1 +\interdelta) \; \sqrt{1+ \delta} \; \sqrt{\frac{2 \Cmax \log
\bigN}{\numobs-\spindex-1} \spec^T (\CovMat_{\Sset\Sset})^{-1} \spec +
\frac{2 \Cmax \sigw \, (1-\frac{\spindex}{\numobs}) \log \bigN}{n
\lambda_{\numobs}^2 }} \\
& \leq & (1+\interdelta) \; \sqrt{1+ \delta} \; \sqrt{\frac{2 \Cmax \;
\spindex \log \bigN}{\numobs-\spindex-1} \frac{1}{\Cmin} + \frac{2
\Cmax \sigw \, \log \bigN}{n \lambda_{\numobs}^2 } }.
\end{eqnarray*}
Now, applying our condition bounding $\numobs, \bigN$ via $\threshbou$
and $\ThreshUp$, we have
\begin{eqnarray*}
\frac{1}{\lambda_{\numobs}} \Exs[\max_j \Zvar_j(\vstar)] & < &
(1+\interdelta) \; \sqrt{1+\delta} \; \sqrt{\matbound^2 \, \left(1 -
\frac{\threshbou \log \bigN}{\numobs-\spindex-1} \right) + \frac{2
\Cmax \sigw \, \log \bigN}{\numobs \lambda_{\numobs}^2 } }.
\end{eqnarray*}
Recall that by assumption, as $\numobs, \bigN \rightarrow +\infty$, we
have that $\frac{\log \bigN}{\numobs \lambda_{\numobs}^2}$ and
$\frac{\log \bigN}{\numobs-\spindex-1}$ converge to zero.
Consequently, the RHS converges to $(1+\interdelta) \,
\sqrt{(1+\delta)} \matbound$ as $\numobs, \bigN \rightarrow \infty$.
Hence, we have
\begin{eqnarray*}
\lim_{\numobs \rightarrow +\infty} \frac{1}{\lambda_{\numobs}}
\Exs[\max_j \Zvar_j(\vstar)] & < & (1+\interdelta) \; \sqrt{1+\delta}
\; \matbound.
\end{eqnarray*}
Since $\interdelta > 0$ and $\delta> 0$ were arbitrary, the result
follows.


\subsection{Proof of Lemma~\ref{LemTwo}}

Consider the function $f: \real^\bigN \rightarrow \real$ given by
\begin{eqnarray*}
f(w) & \defn & \max_{1\leq j \leq \bigN} \left[\sqrt{ \vstar \CondMat}
\; w \right],
\end{eqnarray*}
where $\CondMat \defn \CovMat_{\Sbar \Sbar} - \CovMat_{\Sbar \Sset}
(\CovMat_{\Sset\Sset})^{-1} \CovMat_{\Sset \Sbar}$.  By construction,
for a Gaussian random vector $V \sim N(0,I)$, we have $f(V)
\stackrel{d}{=} \max_{j \in \Sbar} \Zvarcond_j$.

We now bound the Lipschitz constant of $f$.  Let $R =
\sqrt{\CondMat}$. For each $w,v \in \real^\bigN$ and index $j =1,
\ldots, \bigN$, we have
\begin{eqnarray*}
\left |[\sqrt{ \vstar R} w]_j - [\sqrt{ \vstar R}v ]_j \right| & \leq
& \sqrt{\vstar} \left| \sum_{k} R_{jk} [w_k - v_k] \right| \\
& \leq & \sqrt{\vstar} \sqrt{\sum_{k} R_{jk}^2} \; \|w - v\|_2 \\
& \leq & \sqrt{\vstar} \| w - v\|_2,
\end{eqnarray*}
where the last inequality follows since $\sum_{k} R_{jk}^2 =
[\CondMat]_{jj} \leq 1$.  Therefore, by Gaussian concentration of
measure for Lipschitz functions~\cite{Ledoux01,Massart03}, we
conclude that for any $\eta > 0$, it holds that
\begin{eqnarray*}
\Prob[ f(W) \geq \Exs[f(W)] + \eta] & \leq & \exp \left
(-\frac{\eta^2}{2 \vstar} \right), \qquad \mbox{and} \\
\Prob[ f(W) \leq \Exs[f(W)] - \eta] & \leq & \exp \left
(-\frac{\eta^2}{2 \vstar} \right).
\end{eqnarray*}


\subsection{Proof of Lemma~\ref{LemUbehave}}
\label{AppUbehave}
Since the matrix $\Amatt{\Sset}^T \Amatt{\Sset}$ is Wishart with
$\numobs$ degrees of freedom, using properties of the inverse Wishart
distribution, we have $\Exs[(\Amatt{\Sset}^T \Amatt{\Sset})^{-1}] =
\frac {(\CovMat_{\Sset\Sset})^{-1} }{\numobs -\spindex -1}$ (see Lemma
7.7.1 of Anderson~\cite{AndersonStat}).  Thus, we compute
\begin{subequations}
\begin{eqnarray*}
\Exs[\myumean_i] & = &  \frac{-\lambda_\numobs \;
\numobs}{\numobs-\spindex-1} e_i^T \; (\CovMat_{\Sset\Sset})^{-1} \,
\spec, \qquad \mbox{and} \\
\Exs[\myuvar_i] & = & \frac{\sigw}{\numobs} \frac{\numobs}{\numobs -
\spindex-1} \; e_i^T (\CovMat_{\Sset\Sset})^{-1} e_i \; = \;
\frac{\sigw}{\numobs-\spindex-1} e_i^T (\CovMat_{\Sset\Sset})^{-1}
e_i.
\end{eqnarray*}
\end{subequations}
Moreover, using formulae for second moments of inverse Wishart
matrices (see, e.g., Siskind~\cite{Siskind72}), we compute for all
$\numobs > \spindex + 3$
\begin{subequations}
\begin{eqnarray*}
\Exs[\myumean_i^2] & = & \frac{\lambda_\numobs^2 \,
\numobs^2}{(\numobs - \spindex) \, (\numobs - \spindex -3)} \; \left[
\left(e_i^T (\CovMat_{\Sset \Sset})^{-1} \spec\right)^2 +
\frac{1}{\numobs - \spindex -1} \left(\spec^T (\CovMat_{\Sset
\Sset})^{-1} \spec\right) \left(e_i^T (\CovMat_{\Sset \Sset})^{-1}
e_i\right) \right] \\
\Exs[(\myuvar_i)^2] & = & \frac{\sigwfour
\numobs^2}{(\numobs-\spindex-1)^2 \; (\numobs - \spindex) \; (\numobs
- \spindex -3)} \; \left(e_i^T (\CovMat_{\Sset\Sset})^{-1} e_i
\right)^2 \; \left[1 + \frac{1}{\numobs- \spindex-1} \right].
\end{eqnarray*}
\end{subequations}

We now compute and bound the variance of $\myumean_i$.  Setting $A_i =
e_i^T (\CovMat_{\Sset \Sset})^{-1} \spec$ and $B_i = e_i^T
(\CovMat_{\Sset \Sset})^{-1} \spec$ for shorthand, we have
\begin{eqnarray*}
\var(\myumean_i) & = & \frac{\lambda_\numobs^2 \, \numobs^2}{(\numobs
  - \spindex) \, (\numobs - \spindex -3)} \; \left[A_i^2 +
  \frac{1}{\numobs- \spindex -1} A_i B_i \right] -
\frac{\lambda_\numobs^2 \; \numobs^2}{(\numobs-\spindex-1)^2} A_i^2 \\
& = & \frac{\lambda_\numobs^2 \numobs^2}{(\numobs - \spindex) \,
    (\numobs-\spindex-3)} \left[A_i^2 \; \left( 1- \frac{(\numobs-
    \spindex) \, (\numobs - \spindex-3)}{(\numobs - \spindex-1)^2}
    \right) + \frac{1}{\numobs- \spindex -1} A_i B_i \right]  \\
& \leq & 2 \lambda_\numobs^2 \left[ \frac{3 A_i^2}{\numobs- \spindex}
+ \frac{A_i B_i}{\numobs - \spindex -1} \right]
\end{eqnarray*}
for $\numobs$ sufficiently large.  Using the bound $\|(\CovMat_{\Sset
\Sset})^{-1} \|_\infty \leq \Dconmax$, we see that the quantities
$A_i$ and $B_i$ are uniformly bounded for all $i$.  Hence, we conclude
that, for $\numobs$ sufficiently large, the variance is bounded as
\begin{eqnarray}
\label{EqnMyumeanVarBound}
\var(\myumean_i) & \leq & \frac{K \lambda_\numobs^2}{\numobs -
\spindex}
\end{eqnarray}
for some fixed constant $K$ independent of $\spindex$ and $\numobs$.

Now since $|\Exs[\myumean_i]| \leq \frac{2 \Dconmax \lambda_\numobs
\numobs}{\numobs- \spindex-1}$, we have
\begin{equation*}
|\myumean_i - \Exs[\myumean_i]| \; \geq \; |\myumean_i| -
|\Exs[\myumean_i]| \; \geq \; |\myumean_i| - \frac{2 \Dconmax
  \lambda_\numobs \numobs}{\numobs-\spindex-1}.
\end{equation*}
Consequently, making use of Chebyshev's inequality, we have
\begin{eqnarray*}
\Prob[|\myumean_i| \geq \frac{6 \Dconmax \lambda_\numobs
\numobs}{\numobs-\spindex-1}] & = & \Prob[|\myumean_i| - \frac{2
\Dconmax \lambda_\numobs \numobs}{\numobs -\spindex-1}\geq \frac{4
\Dconmax \lambda_\numobs \numobs}{\numobs - \spindex-1}]\\
& \leq & \Prob[|\myumean_i - \Exs[\myumean_i]| \geq \frac{4 \Dconmax
\lambda_\numobs \numobs}{\numobs - \spindex-1}] \\
& \leq & \frac{\var(\myumean_i)}{16 \Dconmax^2 \lambda_\numobs^2} \\
& \leq & \frac{K}{16 \Dconmax \, (\numobs - \spindex)},
\end{eqnarray*}
where the final step uses the bound~\eqref{EqnMyumeanVarBound}.

We now compute and bound the variance of $\myuvar_i$.  We have
\begin{eqnarray*}
\var(\myuvar_i) & = & \frac{\sigwfour
  \numobs^2}{(\numobs-\spindex-1)^2 \; (\numobs - \spindex) \;
  (\numobs - \spindex -3)} \; \left(A_i^2 \; \left[1 +
  \frac{1}{\numobs- \spindex-1} \right] \right) -
\frac{\sigwfour}{(\numobs-\spindex-1)^2} A_i^2 \\
& = & \frac{\sigwfour \numobs^2}{(\numobs-\spindex-1)^2 \; (\numobs -
  \spindex) \; (\numobs - \spindex -3)} \; \left(A_i^2 \; \left[1 +
  \frac{1}{\numobs- \spindex-1} - \frac{(\numobs-\spindex) \; (\numobs
  - \spindex - 3) }{\numobs^2}\right] \right) \\
& \leq & \frac{K \sigwfour}{(\numobs-\spindex-1)^3}
\end{eqnarray*}
for some constant $K$ independent of $\spindex$ and $\numobs$.   Consequently,
applying Chebyshev's inequality, we have
\begin{eqnarray*}
\Prob[\myuvar_i \geq 2 \Exs[\myuvar_i]] \; = \; \Prob[\myuvar_i -
\Exs[\myuvar_i] \geq \Exs[\myuvar_i]] & \leq & 
\frac{\var(\myuvar_i)}{(\Exs[\myuvar_i])^2} \\
& \leq & \frac{K}{(\numobs - \spindex-1)^3}
\frac{1}{\frac{\sigwfour}{\numobs^2} e_i^T (\CovMat_{\Sset
\Sset})^{-1} e_i} \\
& \leq & \frac{K \numobs^2 \Cmax }{\sigwfour (\numobs - \spindex-1)^3} \\
& \leq & \frac{K'}{\numobs - \spindex-1}
\end{eqnarray*}
for some constant $K'$ independent of $\spindex$ and $\numobs$.


\subsection{Proof of Lemma~\ref{LemExpInfinity}}
As in the proof of Lemma~\ref{LemOne}, we define and bound
\begin{eqnarray*}
\Delta_Z(i,j) & \defn & \Exs[(\Zvar_i - \Zvar_j)^2] \; \leq \; 2 \Cmax
\vstar.
\end{eqnarray*}
Now let $(X_1, \ldots, X_\bigN)$ be an i.i.d. zero-mean Gaussian
vector with $\var(X_i) = \Cmax \vstar$, so that $\Delta_X(i,j) \defn
\Exs[(X_i - X_j)^2] = 2 \Cmax \vstar$.  If we set
\begin{eqnarray*}
\Delta^* & \defn & \max_{i,j \in \Sbar} \left |\Delta_X(i,j) -
\Delta_Z(i,j) \right|,
\end{eqnarray*}
then, by applying a known error bound for the Sudakov-Fernique
inequality~\cite{Chatterjee05}, we are guaranteed that
\begin{eqnarray}
\label{EqnKeyLowerBound}
\Exs[\max_{j \in \Sbar} \Zvar_j] & \geq & \Exs [\max_{j \in \Sbar} X_j] -
\sqrt{\Delta^* \log \bigN}.
\end{eqnarray}
We now show that the quantity $\Delta^*$ is upper bounded by
\begin{eqnarray*}
\Delta^* & \leq & 2 \vstar \; (\Cmax - \frac{1}{\Cmax}).
\end{eqnarray*}
Using the inversion formula for block-partitioned
matrices~\cite{Horn85}, we have
\begin{equation*}
\CondMat \; \defn \; \CovMat_{\Sbar \Sbar} - \CovMat_{\Sbar \Sset}
(\CovMat_{\Sset\Sset})^{-1} \CovMat_{\Sset \Sbar} \; = \;
\left[\CovMat^{-1}\right]_{\Sbar \Sbar}.
\end{equation*}
Consequently, we have the lower bound
\begin{eqnarray*}
\Exs[(\Zvar_i - \Zvar_j)^2] & = & \vstar (e_i -e_j)^T \CondMat (e_i -
e_j) \\
& \geq & 2 \vstar \myeigmin(\CondMat) \\
& \geq &  2 \vstar \myeigmin(\CovMat^{-1}) \\
& = & \frac{2 \vstar}{\Cmax}.
\end{eqnarray*}
In turn, this leads to the upper bound
\begin{eqnarray*}
\Delta^* & = & \max_{i, j \in \Sbar} \left | \Delta_X(i,j) -
\Delta_Z(i,j) \right| \\
& = & \max_{i, j \in \Sbar} \left [ 2 \vstar \Cmax - \Delta_Z(i,j) \right] \\
& \leq & 2 \vstar \, \left(\Cmax - \frac{1}{\Cmax} \right).
\end{eqnarray*}

We now analyze the behavior of $\Exs [\max_{j \in \Sbar} X_j]$.  Using
asymptotic results on the extrema of i.i.d. Gaussian
sequences~\cite{GalambosExt,DavNag03}, we have $\lim_{\bigN
\rightarrow +\infty} \frac{\Exs[\max_{j \in \Sbar} X_j]}{\sqrt{ 2
\Cmax \vstar \log \bigN}} = 1$.  Consequently, for all $\interdelta >
0$, there exists an $\bigN(\interdelta)$ such that for all $\bigN \geq
\bigN(\interdelta)$, we have
\begin{eqnarray*}
\Exs[\max_{j \in \Sbar} X_j] & \geq & (1-\interdelta) \sqrt{ 2 \Cmax
\vstar \log \bigN}.
\end{eqnarray*}
Applying this lower bound to the bound~\eqref{EqnKeyLowerBound}, we
have
\begin{eqnarray}
\frac{1}{\lambda_{\numobs}} \Exs[\max_{j \in \Sbar} \Zvar_j] & \geq &
\frac{1}{\lambda_{\numobs}} \left[(1-\interdelta) \, \sqrt{2 \Cmax
\vstar \log \bigN} - \sqrt{ \Delta^* \, \log \bigN} \right] \nonumber
\\
& \geq & \frac{1}{\lambda_{\numobs}} \left[(1-\interdelta) \, \sqrt{2
\Cmax \vstar \log \bigN} - \sqrt{ 2 \, \vstar \, (\Cmax -
\frac{1}{\Cmax}) \, \log \bigN} \right] \nonumber \\
\label{EqnSplit}
& = & \left [(1-\interdelta) \sqrt{\Cmax} - \sqrt{\Cmax -
\frac{1}{\Cmax}} \right] \; \sqrt{ 2
\frac{\vstar}{\lambda_{\numobs}^2} \log \bigN}.
\end{eqnarray}

First, assume that $\lambda_{\numobs}^2/\vstar$ does not diverge to
infinity.  Then, there exists some $\alpha > 0$ such that
$\frac{\lambda_{\numobs}^2}{\vstar} \leq \alpha$ for all sufficiently
large $n$.  In this case, we have from the bound~\eqref{EqnSplit} that
\begin{eqnarray*}
\frac{1}{\lambda_{\numobs}} \Exs[\max_{j \in \Sbar} \Zvar_j] & \geq &
\gamma \sqrt{\log \bigN}
\end{eqnarray*}
where $\gamma \defn \left [(1-\interdelta) \sqrt{\Cmax} - \sqrt{\Cmax
    - \frac{1}{\Cmax}} \right] \; \frac{1}{\sqrt{\alpha}} > 0$.  (Note
that by choosing $\interdelta > 0$ sufficiently small, we can always
guarantee that $\gamma > 0$, since $\Cmax \geq 1$.)  This completes
the proof of condition (b) in the lemma statement.

Otherwise, we may assume that $\lambda_{\numobs}^2/\vstar \rightarrow
+\infty$.  We compute
\begin{eqnarray*}
\frac{1}{\lambda_{\numobs}} \sqrt{ 2 \vstar \log \bigN} & = &
\sqrt{1-\delta} \, \sqrt{\frac{2 \log \bigN}{\numobs-\spindex-1}
\spec^T (\CovMat_{\Sset\Sset})^{-1} \spec + \frac{2 \sigw \,
(1-\frac{s}{\numobs}) \log \bigN}{\numobs \lambda_{\numobs}^2 } } \\
& \geq & \sqrt{1-\delta} \, \sqrt{\frac{2 \log \bigN}{\numobs-\spindex-1} \spec^T
(\CovMat_{\Sset\Sset})^{-1} \spec} \\
& \geq & \sqrt{\frac{1-\delta}{\Cmax}} \, \sqrt{\frac{2 \spindex \log
\bigN}{\numobs- \spindex -1}}.
\end{eqnarray*}
We now apply the condition
\begin{eqnarray*}
\frac{2 \spindex \log \bigN }{\numobs-\spindex-1} & > &
      \frac{1}{\ThreshLow - \threshbou} \; = \; \Cmax \,
      (2-\matbound)^2/\left[ \left[\sqrt{\Cmax} - \sqrt{\Cmax -
      \frac{1}{\Cmax}} \right]^2 -\threshbou \Cmax (2-\matbound)^2
      \right]
\end{eqnarray*}
to obtain that
\begin{eqnarray}
\label{EqnLower}
\frac{1}{\lambda_{\numobs}} \Exs[\max_{j \in \Sbar} \Zvar_j] & \geq &
\sqrt{(1-\delta)} \; \frac{(1-\interdelta) \, \sqrt{\Cmax} - \sqrt{\Cmax-
\frac{1}{\Cmax}}}{\sqrt{\left[\sqrt{\Cmax} - \sqrt{\Cmax -
\frac{1}{\Cmax}}\right]^2 - \threshbou \Cmax \, (2-\matbound)^2}} \;
(2-\matbound)
\end{eqnarray}

Recall that $\threshbou \Cmax \, (2-\matbound)^2 > 0$ is fixed, and
moreover that $\delta, \interdelta > 0$ are arbitrary.  Let $F(\delta,
\interdelta)$ be the lower bound on the RHS~\eqref{EqnLower}.  Note
that $F$ is a continuous function, and moreover that
\begin{eqnarray*}
F(0, 0) & = & \frac{\sqrt{\Cmax} - \sqrt{\Cmax -
\frac{1}{\Cmax}}}{\sqrt{\left[\sqrt{\Cmax} - \sqrt{\Cmax -
\frac{1}{\Cmax}}\right]^2 - \threshbou \Cmax \, (2-\matbound)^2 }} \;
(2-\matbound) \; > \; (2-\matbound).
\end{eqnarray*}
Therefore, by the continuity of $F$, we can choose $\delta,
\interdelta > 0$ sufficiently small to ensure that for some $\gamma >
0$, we have $\frac{1}{\lambda_{\numobs}} \Exs[\max_{j \in \Sbar}
\Zvar_j] \geq (2-\matbound) \; (1 + \gamma)$ for all sufficiently
large $\numobs$.


\subsection{Proof of Lemma~\ref{LemGaussConcenLower}}
\label{AppLemGaussConcenLower}

This claim follows from the proof of Lemma~\ref{LemTwo}.


\bibliographystyle{plain}

\end{document}